\documentclass{article}

\usepackage{PRIMEarxiv}

\usepackage[utf8]{inputenc} % allow utf-8 input
\usepackage[T1]{fontenc}    % use 8-bit T1 fonts
\usepackage{hyperref}       % hyperlinks
\usepackage{url}            % simple URL typesetting
\usepackage{booktabs}       % professional-quality tables
\usepackage{amsfonts}       % blackboard math symbols
\usepackage{amsmath}
\usepackage{nicefrac}       % compact symbols for 1/2, etc.
\usepackage{microtype}      % microtypography
\usepackage{lipsum}
\usepackage{graphicx}
\usepackage{subfigure}
\usepackage{mathpazo}
\usepackage{sectsty}
\usepackage{tikz}

\graphicspath{{media/}}     % organize your images and other figures under media/ folder

\def\bfx{{\mathbf{x}}}
\def\bfy{{\mathbf{y}}}
  
\title{A Scoping Review on Simulation-based Design Optimization in Marine Engineering: Trends, Best Practices, and Gaps
}

\author{
  Andrea Serani$^{1,\star}$, Thomas Scholcz$^2$, and Valentina Vanzi$^3$\\
  $^1$National Research Council-Institute of Marine Engineering, Rome, Italy\\
  $^2$Maritime Research Institute Netherlands, Wageningen, The Netherlands\\
  $^3$Department of Biomedicine and Prevention, University of Rome Tor Vergata, Rome, Italy\\
  $^\star$\texttt{andrea.serani@cnr.it} \\
}

\begin{document}
% \begin{tikzpicture}[remember picture,overlay]
%    \node [rectangle, fill=cyan, fill opacity=0.5, anchor=north, minimum width=\paperwidth, minimum height=3cm, text width=\textwidth, align=center, text height=5ex, text depth=10ex, align=left] at (current page.north) {\sffamily\small 
%    \textbf{This is a preprint of the following article:}\\
%    A. Serani, T. P. Scholcz, V. Vanzi, A Scoping Review on Simulation-based Design Optimization in Marine Engineering: Trends, Best Practices, and Gaps. \textit{Archives of Computational Methods in Engineering}, 2024.\\
%    \textbf{The published article is available by following the DOI: \texttt{10.1007/s11831-024-10127-1}, which may differ from this preprint.}
%    };
% \end{tikzpicture}

\begin{tikzpicture}[remember picture,overlay]
   % Nodo per il riempimento con trasparenza
   \node [rectangle, fill=cyan, fill opacity=0.5, anchor=north, minimum width=\paperwidth, minimum height=3cm] at (current page.north) {};

   % Nodo separato per il testo, senza trasparenza
   \node [anchor=north, minimum width=\paperwidth, minimum height=3cm, text width=\textwidth, align=center, text height=5ex, text depth=10ex, align=left] at (current page.north) {
     \sffamily\small
     \textbf{This is a preprint of the following article:}\\
     A. Serani, T. P. Scholcz, V. Vanzi, A Scoping Review on Simulation-based Design Optimization in Marine Engineering: Trends, Best Practices, and Gaps. \textit{Archives of Computational Methods in Engineering}, 2024.\\
     \textbf{The published article is available by following the DOI: \texttt{10.1007/s11831-024-10127-1}, which may differ from this preprint.}
   };
\end{tikzpicture}

\maketitle

\begin{abstract}
This scoping review assesses the current use of simulation-based design optimization (SBDO) in marine engineering, focusing on identifying research trends, methodologies, and application areas. Analyzing 277 studies from Scopus and Web of Science, the review finds that SBDO is predominantly applied to optimizing marine vessel hulls, including both surface and underwater types, and extends to key components like bows, sterns, propellers, and fins. It also covers marine structures and renewable energy systems. A notable trend is the preference for deterministic single-objective optimization methods, indicating potential growth areas in multi-objective and stochastic approaches. The review points out the necessity of integrating more comprehensive multidisciplinary optimization methods to address the complex challenges in marine environments. Despite the extensive application of SBDO in marine engineering, there remains a need for enhancing the methodologies' efficiency and robustness. This review offers a critical overview of SBDO's role in marine engineering and highlights opportunities for future research to advance the field.
\end{abstract}

% keywords can be removed

\section{Introduction}\label{sec1}
Simulation-based design optimization (SBDO), also known as simulation-driven design optimization {(SDDO)}, has emerged as a critical tool in marine engineering, profoundly impacting various aspects of the field. This approach, which integrates numerical solutions with computer-aided design software and optimization algorithms, empowers engineers to refine performance, cost-efficiency, and safety in marine structures, including ships, underwater vehicles, offshore platforms, and notably, marine energy production systems.

Traditional marine engineering practices, reliant on empirical data and heuristic approaches, often face limitations in adaptability and precision. These methods, though time-tested, struggle to cope with the increasing complexity of marine engineering challenges, especially in the face of stringent environmental regulations and the demand for higher efficiency. SBDO addresses these challenges by enabling a more nuanced exploration of design possibilities, leveraging computational power to identify optimal solutions that balance performance, cost, and environmental considerations.

In ship hull design, SBDO replaces traditional methods, which are heavily reliant on experience and trial-and-error approaches. By analyzing hydrodynamic performance across different hull designs, SBDO enables the optimization of shape and dimensions, thus reducing drag and enhancing fuel efficiency \cite{lowe_automatic_1994, campana_shape_2006, peri_design_2001}.

For marine propulsion systems, SBDO is invaluable in dealing with the complexity of various components like engines, propellers, shafts, and rudders. It facilitates the optimization of these components for maximum efficiency and reduced fuel consumption \cite{vesting_surrogate_2014, ma_design_2014, chen_parametric_2018, nouri_optimization_2018, mirjalili_confidence-based_2018, diez_hydroelastic_2012, favacho_contribution_2016, esmailian_systematic_2017, lu_research_2021}.

A pivotal area where SBDO is making significant strides is in the development and optimization of marine energy production systems. As the world increasingly seeks sustainable energy sources, marine energy systems, such as tidal \cite{kinnas_computational_2012, zhang_optimization_2016, shi_optimal_2015, huang_optimization_2016, sun_prediction_2019, im_duct_2020, khanjanpour_optimization_2020,ambarita_computational_2021, yeo_tidal_2022} and wave energy converters \cite{e_silva_hydrodynamic_2016, simonetti_optimization_2017, tao_optimized_2021, bao_parametric_2021}, have gained prominence. SBDO plays a crucial role in designing these systems to maximize energy extraction and efficiency while ensuring resilience to marine environmental challenges. The optimization of these systems is vital for advancing renewable energy technologies and contributes significantly to sustainable marine practices.

Additionally, SBDO enhances the safety and reliability of marine structures. For offshore structures \cite{yang_robust_2015}, which face harsh environmental conditions, SBDO is instrumental in evaluating and improving structural integrity under various scenarios.

%Overall, SBDO is driving innovation in marine engineering, improving design quality, accelerating development timelines, and reducing costs. This method's ability to identify and solve potential design issues early on is revolutionizing the field. This scoping review aims to present a comprehensive, current overview of SBDO in marine engineering, highlighting its applications, including marine energy systems, and pointing future research directions within marine and ocean engineering contexts.

Looking ahead, the field of SBDO in marine engineering is poised for significant advancements. Emerging trends like the integration of machine learning algorithms and the incorporation of real-time data analytics are expected to further revolutionize SBDO applications. These advancements will not only refine the optimization process but also open new avenues for addressing complex, multifaceted marine engineering challenges. This scoping review aims to present a comprehensive, current overview of SBDO in marine engineering, highlighting its applications and pointing to future research directions within marine and ocean engineering contexts.

%=================================================================================
\section{Scoping Review Methodology}
Due to a noticeable increase in research output and the proliferation of primary research over the past few years, the need to systematically identify and synthesize the existing literature has become mandatory in research. This critical issue has first arisen in clinical medicine but nowadays it represents a priority in many other disciplines including engineering \cite{brereton2007lessons}. Scoping reviews are extremely useful to accomplish this goal. The original framework for conducting scoping reviews was proposed by Arksey and O’Malley in 2005 \cite{arksey2005scoping} and further extended by Joanna Briggs Institute (JBI) Collaboration in 2015 \cite{peters2015guidance}. Recently, the JBI Scoping Reviews Methodology Group formally defined scoping reviews as a ``\emph{type of evidence synthesis that aims to systematically identify and map the breadth of evidence available on a particular topic, field, concept, or issue, often irrespective of source (i.e., primary research, reviews, non-empirical evidence) within or across particular contexts}'' \cite{munn2022scoping}. 
Despite other review methods, scoping reviews use a broader approach for mapping literature and addressing a broader research question without performing articles' quality assessment \cite{sharma2023write}.

\subsection{Research Questions}
Central to this review is the exploration of current best practices in SBDO applied to marine engineering. This inquiry is structured into three fundamental questions:
\begin{enumerate}
\item What are the primary aims and approaches in the existing literature on SBDO methods in marine engineering, and how do they compare?
\item What issues are encountered when applying SBDO methods to marine engineering problems?
\item What are the main research gaps and potential future directions in this field?
\end{enumerate}
%-------------------------------------------------------------
\begin{figure*}[!b]
    \centering
    \includegraphics[width=1\textwidth]{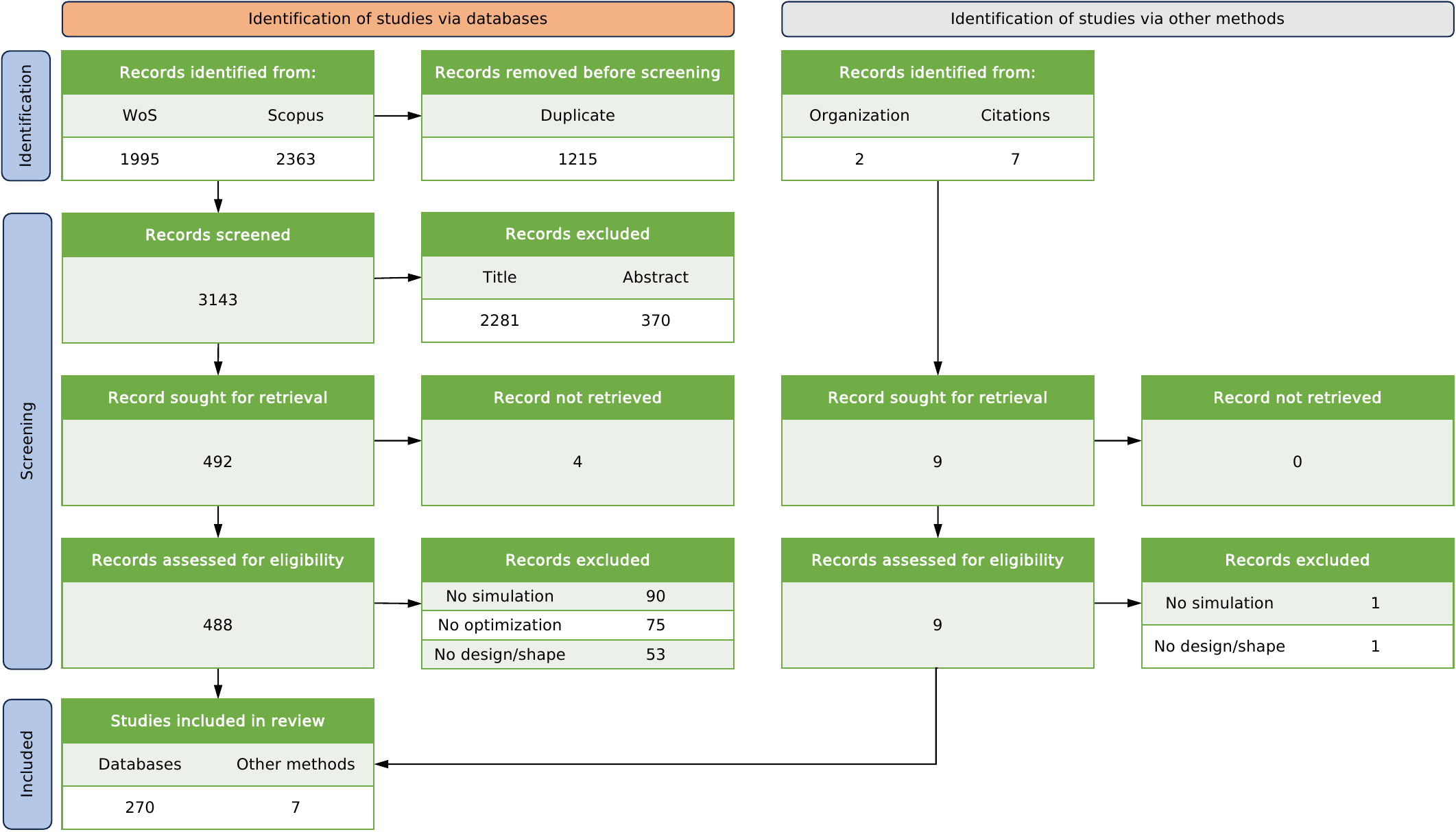}
    \caption{PRISMA flow chart.}
    \label{fig:prisma}
\end{figure*}

\subsection{Inclusion and Exclusion Criteria}
The inclusion criteria for the articles in this review were meticulously defined to ensure a focused and relevant collection of literature. Articles were selected based on their direct relevance to SBDO applications in marine engineering. This included studies demonstrating the use of SBDO in practical marine engineering projects, theoretical advancements in SBDO methods specific to marine applications, and reviews of SBDO methods within the marine engineering context.

Exclusion criteria were equally stringent to maintain the review's scope and quality. Articles not directly related to SBDO, such as those focusing on general design optimization without a clear simulation-based component, were excluded. Studies outside the realm of marine engineering, or those employing SBDO in a manner not applicable to marine engineering challenges, were also omitted. 
Furthermore, non-peer-reviewed articles, such as conference abstracts/papers and editorials, were excluded to ensure the review's academic rigor.

%\subsection{Databases and keywords}
%This study used Web of Science (WoS) and Scopus as databases, and a number of grouped keywords to cover the full extent of current research in the field. The following were the search terms employed: (''Simulation*'' OR ''Computation*'') AND (''Optimi*’'') AND (''Design*'' OR ''Shape*'' OR ''Form*'') AND (''Ship*'' OR ''Hull'' OR ''Vessel'' OR ''Marine'' OR ''Ocean'').

\subsection{Databases and Keywords}
Web of Science (WoS) and Scopus were chosen as the primary databases for their extensive coverage of interdisciplinary scientific literature, ensuring a comprehensive collection of relevant articles in marine engineering and optimization. These databases are renowned for their rigorous indexing of high-quality, peer-reviewed academic journals, which aligns with the review's emphasis on academic rigor. %While other databases like Google Scholar and IEEE Xplore were considered, they were ultimately excluded due to their broader, less specialized focus, which could potentially dilute the review's specificity in SBDO within marine engineering.

%The search strategy employed grouped keywords to encompass the full spectrum of SBDO research in marine engineering. The keywords used were: (``Simulation*'' OR ``Computation*'') AND (``Optimi*'') AND (``Design*'' OR ``Shape*'' OR ``Form*'') AND (``Ship*'' OR ``Hull'' OR ``Vessel'' OR ``Marine'' OR ``Ocean''). 
{The bibliographic search strategy was carefully designed to capture the broad scope of SBDO research in marine engineering, employing a combination of keywords specifically targeted within the titles, abstracts, and keywords sections (\texttt{TITLE-ABS-KEY}) of articles. The chosen keywords aimed to include a comprehensive range of studies relevant to the field: \texttt{(``Simulation*'' OR ``Computation*'') AND (``Optimi*'') AND (``Design*'' OR ``Shape*'' OR ``Form*'') AND (``Ship*'' OR ``Hull'' OR ``Vessel'' OR ``Marine'' OR ``Ocean'')}. This strategic choice ensured the inclusion of pertinent research while maintaining a focused scope on SBDO applications within marine engineering.}

%----------------------------------------------------------------
\begin{figure*}[!b]
    \centering
    \subfigure[]{\includegraphics[width=0.495\textwidth]{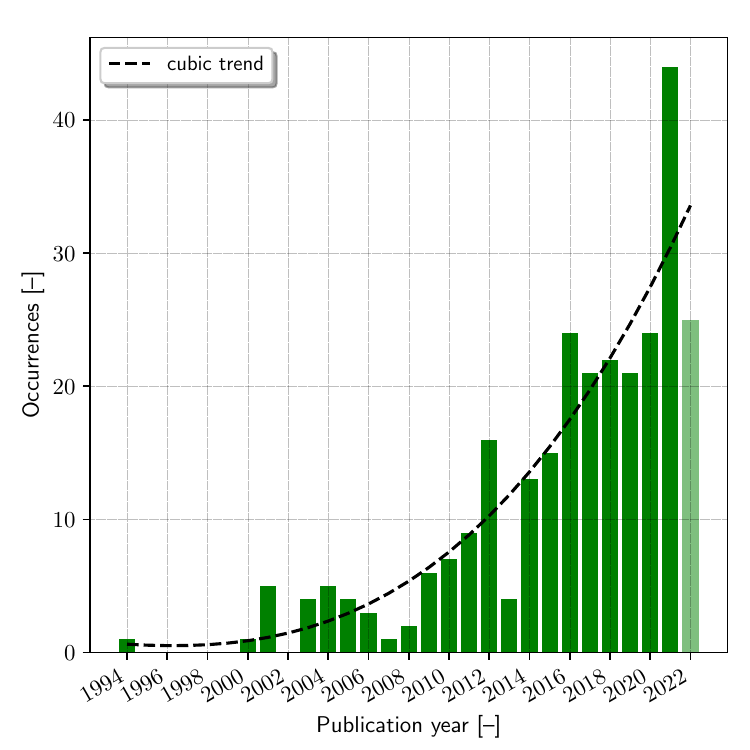}}
    \subfigure[]{\includegraphics[width=0.495\textwidth]{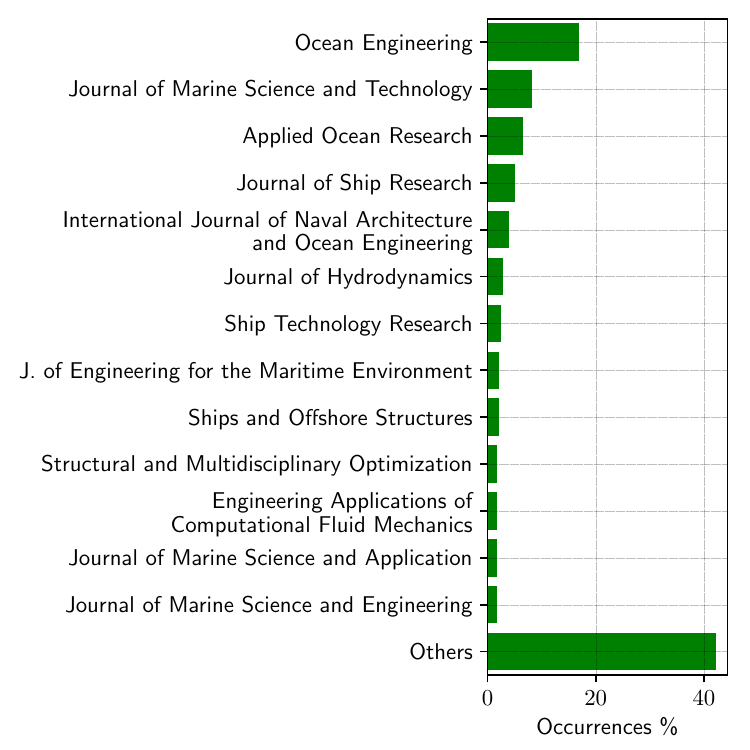}}
    \caption{Publication trend (a) and journals occurrence (b).}
    \label{fig:years}
\end{figure*}
\subsection{Search Procedure}
The preferred reporting items for systematic reviews statement extended to scoping reviews (PRISMA-ScR) are used as reporting guidelines \cite{tricco2018prisma}.
%This paper used the preferred reporting items for systematic reviews statement extended to scoping reviews (PRISMA-ScR) as reporting guidelines \cite{tricco2018prisma}. 
The PRISMA flow diagram (see Fig. \ref{fig:prisma}) meticulously outlines the process undertaken for the selection of articles in the present scoping review. The articles search was conducted on August 1st, 2022, with no restriction on the date of publication and type of study, but considering only journal papers written in English.
The diagram begins with the identification phase, where 3143 records were sourced through WoS and Scopus, indicating a comprehensive initial search strategy. Reference lists of all included articles were scanned to look for literature that had not been obtained previously.

Subsequent stages in the diagram reflect the screening and eligibility assessment processes. Notably, a significant number of records were excluded during the initial screening, likely due to title (2281) and abstract (370) relevance checks. This highlights the precision of our inclusion criteria, ensuring that only the most pertinent articles were considered (492) for full-text review.

The eligibility phase, as depicted, involved a more detailed review of the full texts, leading to further exclusion of articles that did not meet the specific criteria set for this review. These criteria were crucial in filtering out articles that did not include simulation, optimization strategies, or design/shape optimization.

Finally, the included studies (277), as shown in the diagram, represent a curated collection of articles that passed through this rigorous selection process, ensuring a high degree of relevance and quality in the research articles selected for this review.

%=================================================================================
\section{Results}
%This section presents the results of the scoping review. The analysis is structured to provide a holistic understanding of the field's progression, encapsulating trends in publications, the distribution of research across journals, and a deep dive into the core components of SBDO methodologies.
The following subsections delineate the comprehensive findings of the scoping review, focusing on the key developments and trends within the realm of SBDO in marine engineering. This analysis aims to distill a broad spectrum of research efforts into discernible patterns, offering insights into the evolution, current practices, and future directions in the field. By examining a variety of aspects, from publication trends and journal distributions to the nuanced details of optimization techniques and application areas, this section endeavors to provide a holistic understanding of the state-of-the-art in SBDO as applied to marine engineering.

{It may be noted that different terms have been used interchangeably to describe the overarching process of integrating computational simulations with design optimization in marine engineering. While SBDO and SDDO are prevalent, the analysis reveals both their widespread use and nuanced differences.
SBDO emerges as the most comprehensive term, encompassing the full spectrum of leveraging simulation tools for optimizing design parameters. This terminology aligns with the holistic approach of using simulations to inform and drive the optimization process, where the objective is to enhance design performance metrics while navigating through the constraints imposed by complex marine engineering challenges.
On the other hand, SDDO often highlights the initial stages of the design process, where simulations guide the conceptual and preliminary design decisions before formal optimization techniques are applied. This term underscores the importance of simulations in shaping the design space and influencing early design choices, which are crucial for setting the stage for subsequent optimization.
The review suggests that while these terms broadly address the same domain of integrating simulations with optimization, they can reflect different focuses or stages within the broader SBDO process. This distinction is vital for understanding the scope and emphasis of various studies within the field, as well as for appreciating the multifaceted nature of SBDO in marine engineering.}
%----------------------------------------------------------------
\begin{figure*}[!t]
    \centering
    \includegraphics[width=1\textwidth]{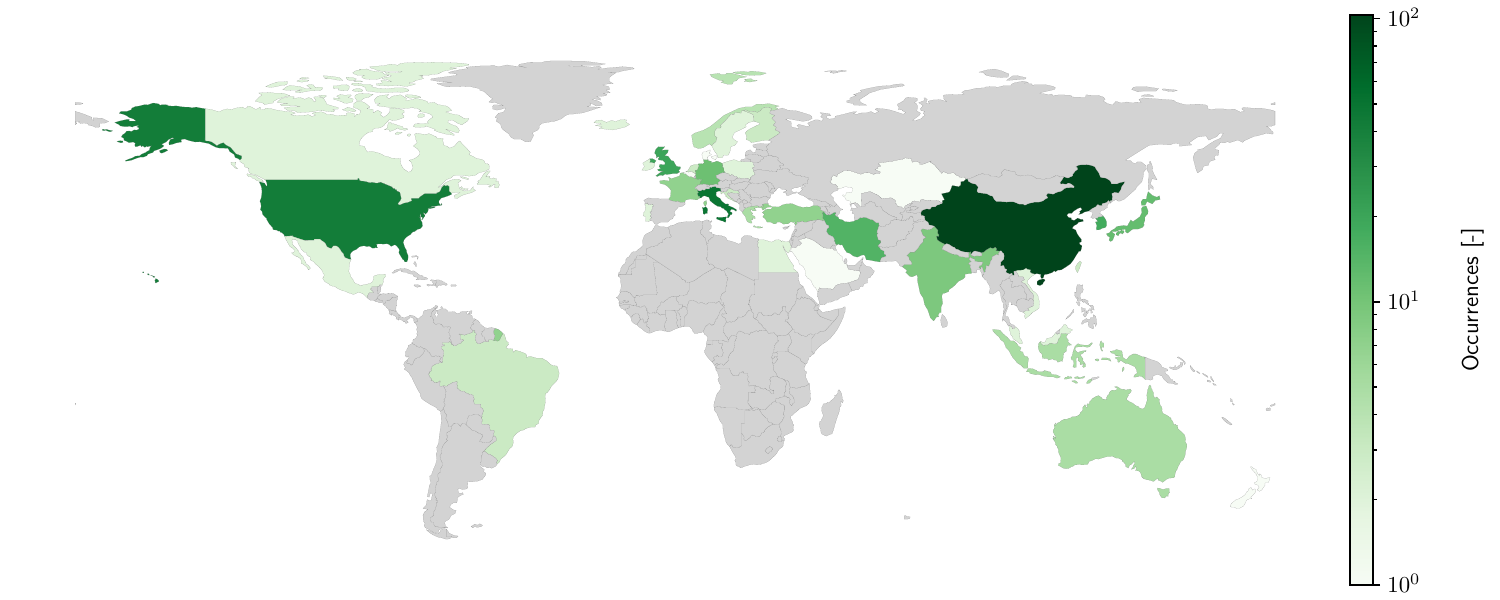}
    \caption{{Publications occurrences geographical distribution (absolute value per country on a logarithmic scale).}}
    \label{fig:world}
\end{figure*}

Figure \ref{fig:years}a illustrates a chronological trend in the number of publications per year on the topic. Starting from 1994, the year of the first publication retrieved on the topic \cite{lowe_automatic_1994}, a noticeable increase in publications can be observed over the years (specifically starting from 2009), indicating a growing interest and advancement in the field. It's important to note that the data for the year 2022 is partial, as the bibliographic research was conducted on August 1, 2022. This uptick reflects the evolving complexity and significance of SBDO in addressing contemporary challenges in marine engineering. The progressive increase underscores the technology’s rising relevance, potentially correlating with advancements in computational capabilities and the growing demand for efficient, optimized marine systems.

\begin{figure}[!b]
    \centering
    \includegraphics[width=0.5\columnwidth]{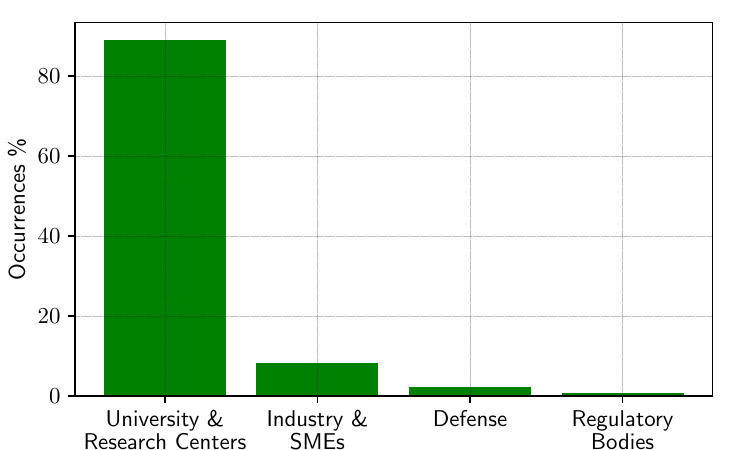}
    \caption{{Publications origin occurrences by entity.}}
    \label{fig:origin}
\end{figure}

Figure \ref{fig:years}b presents a distribution of publications across various journals, highlighting those with the highest frequency of articles. Overall the Ocean Engineering journal covers 17.2\% of the overall publications, whereas  the other journals all contain less than 10\% of the publications on SBDO. Moreover, the category 'Others' encapsulates a range of journals that individually contribute to less than 2\% of total publications, signifying a wide dissemination of research in this field across diverse scientific platforms. This distribution not only reflects the interdisciplinary nature of the field but also points to the key academic outlets that are central to the dissemination of SBDO research.

{Based on a detailed analysis of the distribution and contributions, the results offer intriguing insights into global research trends and collaborative dynamics. The geographical distribution (see Fig. \ref{fig:world}) showcases a significant concentration of contributions from China, accounting for 29.3\% of the papers reviewed, with a diverse representation from 48 different entities. This is followed by Italy (13.9\%), the United States (11.9\%), the United Kingdom (5.7\%), South Korea (5.1\%), Iran (4.3\%), Japan (3.4\%), and Germany (3.1\%), highlighting a global interest and varied focus across these regions. The predominance of university and research center contributions, with 89\% of the instances (see Fig. \ref{fig:origin}), signifies the academic inclination of SBDO research, whereas the industry and small and medium enterprises (SMEs), defense agencies, and regulatory bodies' engagement, though lesser in number, underscore the multi-sectoral relevance of SBDO applications in marine engineering.
This diverse geographical and institutional representation underscores the universal appeal and applicability of SBDO techniques across different marine engineering challenges, reflecting a rich picture of research efforts aimed at advancing marine technology and sustainability. The data suggest a vibrant and collaborative research ecosystem, with significant contributions emerging from both academia and industry, pointing towards an integrated approach to innovation in marine engineering through SBDO.}

The following subsections present a categorization of SBDO research into several key areas, resulting in a systematic description of the vast body of work in this domain. The examination begins with problem formulation strategies, identifying the complex nature and challenges of the design optimizations present in the various studies. Subsequent analysis delves into the parameterization techniques used in SBDO. The focus then shifts to the solvers utilized in SBDO and optimization strategies. Finally, a deeper discussion of the applications is given.

\begin{figure*}[!b]
    \centering
    \subfigure[]{\includegraphics[width=1\textwidth]{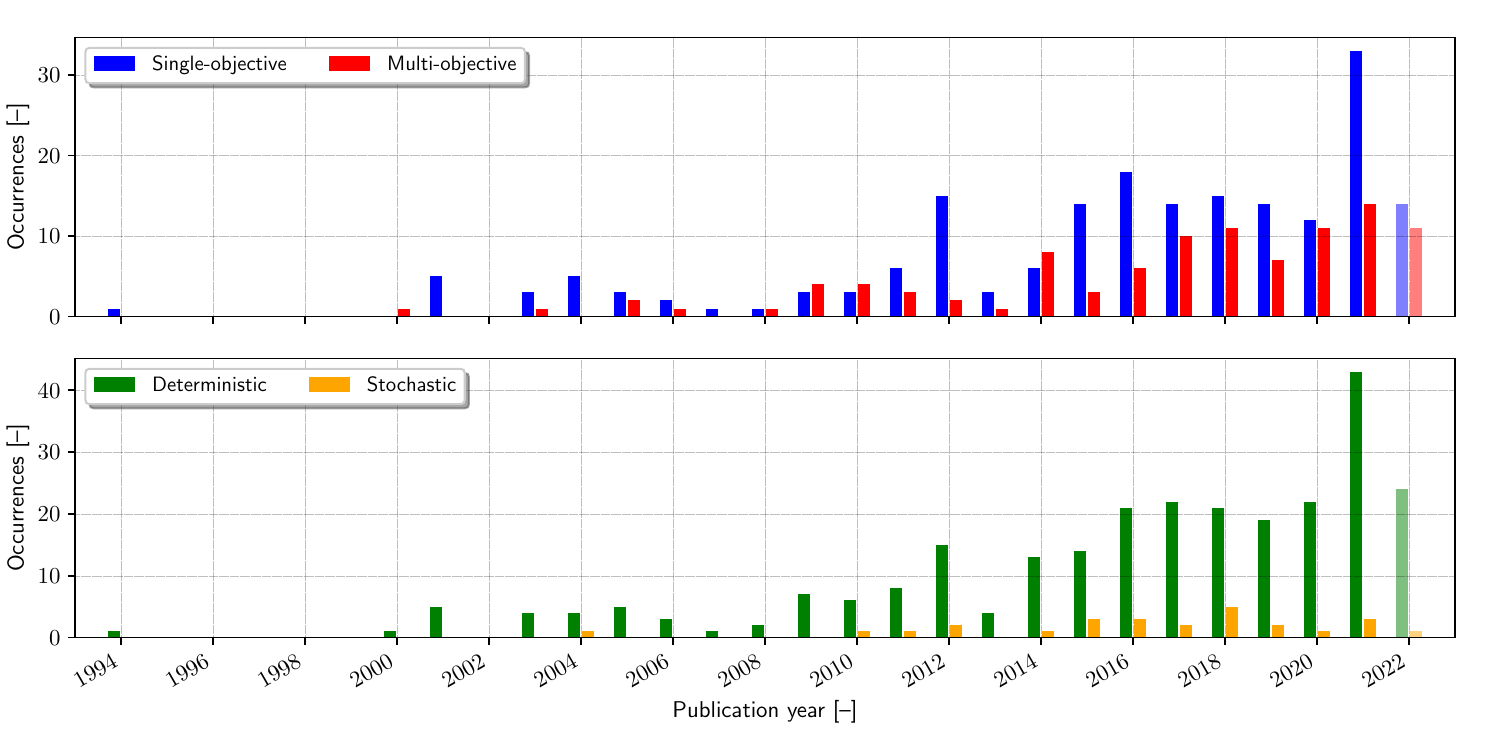}}\\
    \subfigure[]{\includegraphics[width=0.495\textwidth]{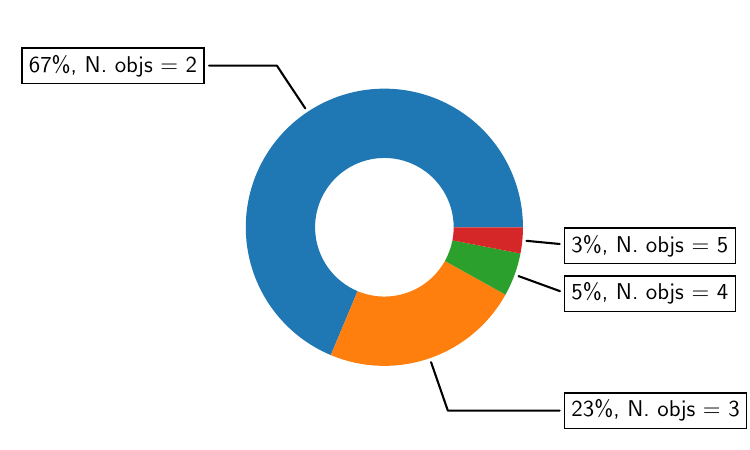}}
    \subfigure[]{\includegraphics[width=0.495\textwidth]{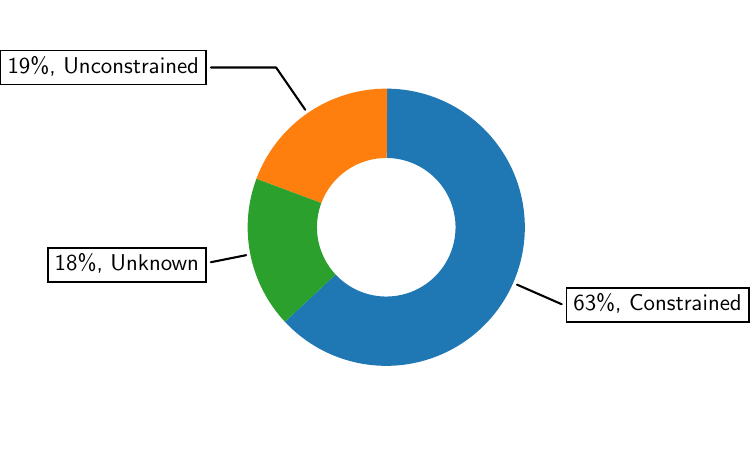}}
    \caption{Problem formulation: (a) occurrences by year (top) single- versus multi-objective and (bottom) deterministic versus stochastic; (b) number of objectives overall occurrence for multi-objective problems; (c) use of constraints overall occurrence.}
    \label{fig:formulations}
\end{figure*}
\begin{figure*}[!b]
    \centering
    \includegraphics[width=0.8\textwidth]{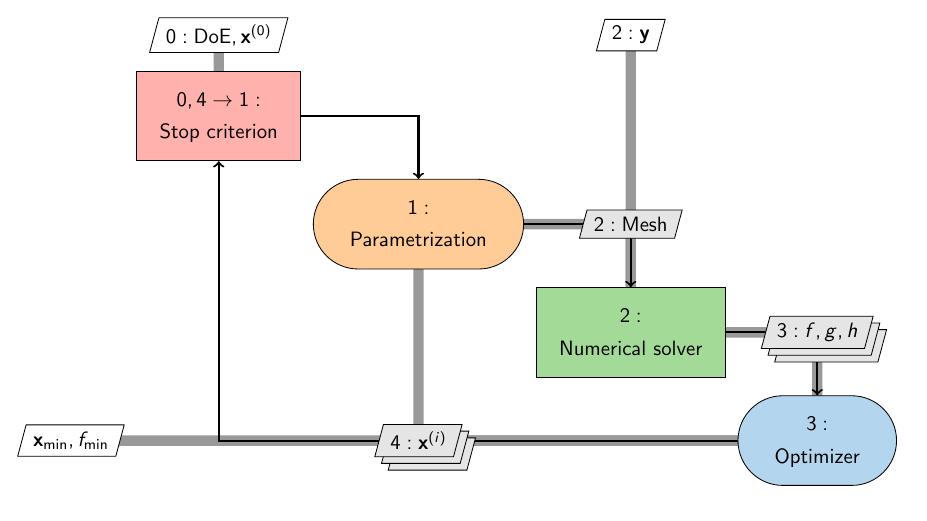}
    \caption{Overview of a general SBDO process through the XDSM diagram.}
    \label{fig:xdsm}
\end{figure*}
%=================================================================================
\subsection{Problem Formulations}

The field of SBDO in marine engineering exhibits a range of problem formulations, from straightforward deterministic single-objective optimization to more complex multi-objective and stochastic optimization approaches. The evolution towards embracing these complexities is gradual, reflecting a preference for simpler, more intuitive methods (see Fig. \ref{fig:formulations}).

Central to the SBDO approach is the deterministic single-objective optimization, which remains predominant due to its clear and straightforward formulation:

\begin{align}\label{eq:soform}
\min_{\bfx} \quad & f(\bfx,\bfy) \\
\text{subject to} \quad & g_i(\bfx,\bfy) \leq 0, \quad i = 1, \ldots, m \nonumber \\
\text{and to} \quad & h_j(\bfx,\bfy) = 0, \quad j = 1, \ldots, p \nonumber\\
\text{and to} \quad & \bfx_l\leq\bfx\leq\bfx_u . \nonumber
\end{align}    
This formulation, with $f$ as the objective function, $\bfx$ as the design variables {(with $\bfx_l$ and $\bfx_u$ the lower and upper bounds)}, $\bfy$ as the environmental and/or operational conditions, $g_i$ as inequality constraints, and $h_j$ as equality constraints, is favored for its ability to produce clear and concise results, making it highly suitable for demonstrating new SBDO methodologies in marine engineering.

Despite the potential to address a broader spectrum of design criteria, the uptake of multi-objective optimization, that reformulate the problem in Eq. \ref{eq:soform} as follows 
\begin{align}\label{eq:moform}
\min_{\bfx} \quad & \{f_1(\bfx,\bfy), f_2(\bfx,\bfy), \ldots, f_k(\bfx,\bfy)\} \\
\text{subject to} \quad & g_i(\bfx,\bfy) \leq 0, \quad i = 1, \ldots, m  \nonumber\\
\text{and to} \quad & h_j(\bfx,\bfy) = 0, \quad j = 1, \ldots, p \nonumber\\
\text{and to} \quad & \bfx_l\leq\bfx\leq\bfx_u, \nonumber
\end{align}
is cautious (see Fig. \ref{fig:formulations}a, top).
This approach, involving the simultaneous optimization of multiple conflicting $k$ objectives, faces challenges due to its increase in required computational resources and complexity.
Figure \ref{fig:xdsm} presents a comprehensive depiction of the SBDO process using the extended design structure matrix (XDSM) \cite{lambe2012extensions}. This representation includes the three main blocks (shape parametrization, numerical solver, and optimizer) of the process, including also a stopping criteria, which may encompass either the convergence of the optimization method or constraints imposed by a limited computational budget.

The adoption of stochastic optimization (see Fig. \ref{fig:formulations}a, bottom), which factors in uncertainty and variability, is still limited. Techniques like robust design optimization (RDO) \cite{yang_robust_2015,lee_surrogate_2016,mirjalili_confidence-based_2018,yin_numerical_2019,bacciaglia_controllable_2021,du_identification_2022}, that focus on performance stability under uncertainty, reliability-based design optimization (RBDO) \cite{young_reliability-based_2010,yang_metamodel_2011,choi_sampling-based_2015,hou_hull_2018,tao_optimized_2021}, which emphasizes safety and reliability standards under probabilistic uncertainty models, and reliability-based and robust design optimization (RBRDO) \cite{pellegrini_formulation_2017,diez_stochastic_2018,serani_hull-form_2022}, that combines RDO and RBDO approaches to ensure that a design is both robust against variability and reliable in terms of meeting safety or success criteria, are not yet widespread, pointing to a significant potential area of development in the field, representing only 9\% of the existing literature.

Figure \ref{fig:formulations}a clearly illustrates the continued preference for single-objective over multi-objective optimization (top) and deterministic over stochastic optimization (bottom) approaches in the marine engineering domain. These preferences underscore the field's inclination towards methodologies that offer straightforward applicability and simplicity. Figure \ref{fig:formulations}b, on the other hand, reveals a modest but growing interest in multi-objective optimization, with a limit to the number of objectives, indicating a cautious approach to embrace complexity in optimization challenges. Examples of many-objectives optimization problems (number of objectives greater than 3) are given in \cite{he_multidisciplinary_2011,kamarlouei_multi-objective_2014,lin_hull_2018,wang_aerodynamic_2021,doijode_machine_2022-1} for 4 objectives and in \cite{peri_high-fidelity_2005,lu_hydrodynamic_2019,mittendorf_hydrodynamic_2021} for 5 objectives.

Furthermore, the analysis of problem formulations in SBDO studies, as depicted in Fig. \ref{fig:formulations}c, reveals that a significant majority of problems (63\%) are formulated with constraints. This indicates that complex real-world conditions and requirements are typically encountered in marine engineering applications. Constraints in SBDO may originate from design, regulatory and safety requirements, physical limitations, and environmental considerations. 

%(MAYBE ADD SOME EXAMPLE REFERENCES ABOUT THE KIND OF CONSTRAINTS...)

The predominance of constrained problems underscores the need for optimization methodologies that can effectively account for these limitations, balancing the achievement of design objectives with adherence to constraint boundaries. Interestingly, a notable 19\% of the problems are identified as unconstrained. This suggests scenarios where design freedom is less restricted, possibly in more theoretical or exploratory studies, or in cases where the primary focus is on optimizing a single aspect of design without the need for balancing it against other factors. Another possibility is the use of implicit geometrical constraints, such that they don't need to be considered in the problem formulation anymore because they are satisfied by definition.   
However, Fig. \ref{fig:formulations}c also highlights a critical gap in current SBDO research - a lack of clarity or information regarding the problem formulation in 18\% of the papers. This ambiguity in the formulation, specifically the absence of clear statements on whether the problems are constrained or not, points to a potential oversight in the documentation or conceptualization of SBDO studies. It raises questions about the comprehensiveness and depth of problem understanding in these cases. The absence of explicit mention of constraints may lead to challenges in replicating or building upon the research, as the constraints (or lack thereof) significantly influence the optimization process and outcomes. 
Furthermore, the figure brings to light an important aspect of SBDO that appears to be insufficiently addressed: the strategies for dealing with constraints. Effective constraint handling is crucial in SBDO, as it directly impacts the feasibility and practicality of the optimized solutions. The lack of detailed discussion on constraint management techniques in a considerable number of studies suggests a need for more focused research in this area. This includes the development and application of advanced constraint-handling techniques, which are essential for ensuring that the solutions generated by SBDO are not only optimal in a mathematical sense but also viable and effective in real-world applications.

%(NEED SOME REFERENCES DESCRIBING CONSTRAINTS HANDLING)
%\textcolor{blue}{FOUND THIS: 
%- A Review on Constraint Handling Techniques for Population‑based Algorithms: from single‑objective to multi‑objective optimization
%- Constraint Handling in Bayesian Optimization -- A Comparative Study of Support Vector Machine, Augmented Lagrangian and Expected Feasible Improvement} 

%\subsubsection{Multidisciplinary design optimization}
The scoping review has finally highlighted a notably sparse yet significant application of multidisciplinary design optimization (MDO) methodologies within the broader context of SBDO in marine engineering, encompassing only about 8\% of the studies. This is particularly noteworthy in a field inherently requiring integration across various disciplines such as hydrodynamics, structural engineering, and materials science for optimal design solutions.
MDO problems focusing on resistance/powering and seakeeping performance improvement have been addressed in the context of various vessels, including surface combatant \cite{peri_multidisciplinary_2003,serani_hull-form_2022},  frigate \cite{peri_multidisciplinary_2005}, and multi-hulls \cite{diez_stochastic_2018,nazemian_multi-objective_2021}. These studies highlight the application of MDO in enhancing specific performance parameters of marine vehicles.
A multilevel hierarchy system approach, which allows for the integration of results from synthesis-level optimization into subsystem optimization and overall coordination of multi-level design systems, was demonstrated in studies like \cite{besnard_constructive_2007} and \cite{hefazi_multidisciplinary_2010}. These works employed methods like constructive artificial neural networks for the MDO of twin H-body vessels and multi-hulls, considering objectives and constraints related to cavitation, structural integrity, stability, hull forms, weights, costs, and payload capacity.
System-level MDO, considering seakeeping, maneuvering, and resistance assessment, was explored in \cite{he_multidisciplinary_2011}, showcasing the comprehensive nature of this MDO approach. In contrast, a generalized collaborative optimization (CO) method for resistance optimization of small water-plane area twin hull (SWATH) vessels was proposed in \cite{xiao_generalised_2012}, signifying the adaptability of CO in focused optimization tasks.
The optimization of an autonomous underwater vehicle (AUV) for various performance metrics such as rapidity, maneuverability, resistance, and energy consumption through CO was undertaken in studies like \cite{luo_application_2015} and \cite{luo_hull_2021}. Additionally, a modified bi-level integrated system collaborative optimization for resistance and weight reduction of a SWATH was proposed in \cite{jiang_modified_2016}.
The application of a multi-objective MDO based on the all-at-once architecture for weight minimization and endurance maximization of an AUV was demonstrated in \cite{liu_multiple_2017}. Resistance optimization and wake flow uniformity of an offshore aquaculture vessel were addressed in \cite{feng_multidisciplinary_2018}, while \cite{zhang_optimum_2019} utilized a concurrent subspace design method for comprehensive MDO of an AUV, covering hull form, structure, propulsion, energy, maneuverability, and general arrangement.
Further studies explored a range of MDO applications \cite{seth_amphibious_2020}, from hydrostructural optimization \cite{young_reliability-based_2010,diez_hydroelastic_2012,garg_high-fidelity_2017} to energy consumption minimization \cite{chen_gradient-based_2018}, showcasing the diversity of MDO applications in marine engineering, employing various architectural approaches such as fluid-structure interaction coupling \cite{diez_experimental_2022}, super element-based multi-level analysis \cite{sun_new_2011}, and uncertainty quantification in system-level MDO \cite{leotardi_variable-accuracy_2016}.

\begin{figure*}[!b]
    \centering
    \subfigure[]{\includegraphics[width=0.495\textwidth]{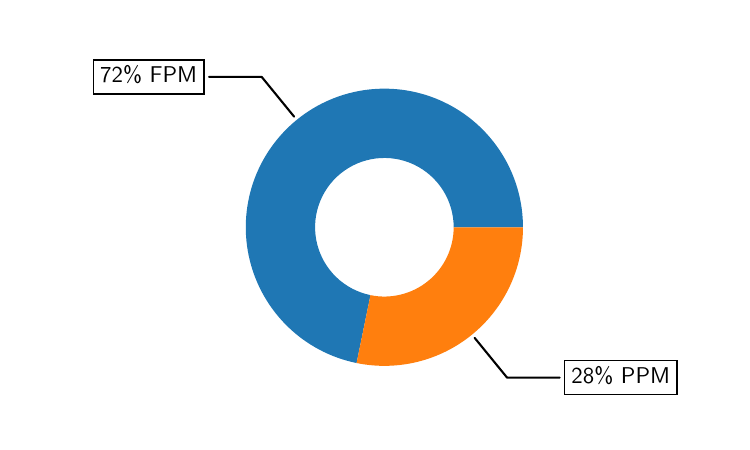}}
    \subfigure[]{\includegraphics[width=0.495\textwidth]{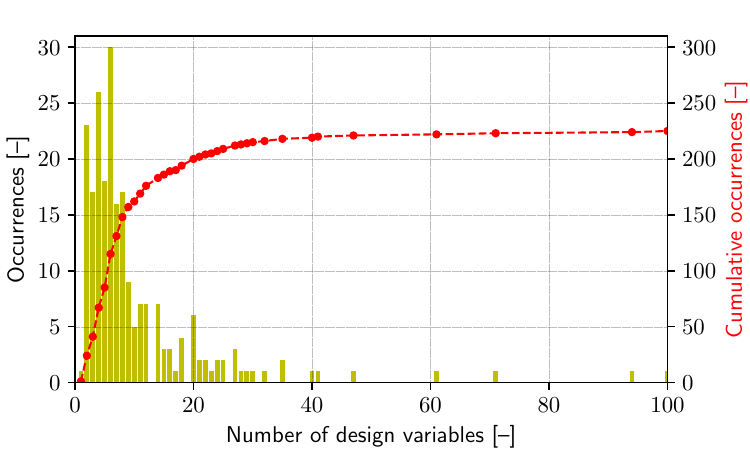}}
    \caption{Occurrences of (a) fully- versus partially-parametric modeling for shape modification and (b) distribution of design-space dimensionality.}
    \label{fig:design_param}
\end{figure*}
%

%=================================================================================
\subsection{Design-space Parameterization}

In the realm of SBDO, the parametrization of the design space is a critical step that significantly influences the optimization process. Parametrization can be categorized broadly into fully-parametric (FPM) and partially-parametric models (PPM) \cite{harries2019faster}. FPMs define every aspect of the design using parameters, offering high control and predictability. PPMs, however, combine parametric elements with non-parametric or fixed aspects, providing a balance between control and flexibility. This distinction is crucial in SBDO, where the choice of parametrization technique impacts the feasibility, efficiency, and scope of the optimization task.

Figure \ref{fig:design_param}a shows the predominant preference for FPM, accounting for 72\%.
This dominance suggests a trend towards well-defined, controlled, and interpretable approaches in design variable specification. FPM approaches include CAD-based \cite{vasudev_modular_2016}, analytical \cite{bagheri_genetic_2014,bagheri_optimizing_2014}, scaling \cite{hefazi_multidisciplinary_2010}, sectional area curves \cite{park_hull-form_2015,kim_hull-form_2016,park_hull_2022}, partial differential equations \cite{lowe_automatic_1994}, Ferguson \cite{choi_sampling-based_2015}, Legendre \cite{coiro_diffuser_2016}, Bezier curves \cite{yang_hydrofoil_2012,koziel_simulation-driven_2012,luo_multi-point_2014,chrismianto_development_2015,sun_prediction_2019,yin_numerical_2019} and surfaces \cite{peri_design_2001,leotardi_variable-accuracy_2016}, Splines \cite{dejhalla_numerical_2001,suzuki_studies_2005}, B-splines \cite{chen_inverse_2004,grigoropoulos_hull-form_2010,mohamad_ayob_uncovering_2011,bertetta_cpp_2012,guha_amitava_application_2015,gaggero_design_2017,nouri_optimization_2018,bonfiglio_improving_2018,lu_hydrodynamic_2019,gaggero_reduced_2020,furcas_design_2020,maia_computational_2021}, T-splines \cite{kostas_ship-hull_2015}, F-splines \cite{liu_cross-entropy_2021}, NURBS \cite{percival_hydrodynamic_2001,zakerdoost_evolutionary_2013,luo_design_2017,tomasz_abramowski_energy_2017,chen_parametric_2018,barbaric_investigation_2020}, PARSEC \cite{ma_hard_2019}, Lackeby \cite{sarioz_inverse_2009,guo_cfd-based_2020}, and Akima \cite{pehlivan_solak_multi-dimensional_2020}. On the other hand, PPM methods such as free-form deformation (FFD) \cite{duvigneau_role_2003,campana_numerical_2009,tahara_single-_2011,li_multiobjective_2013,vasile_ship_2013,li_bow_2014,diez_design-space_2015,garg_high-fidelity_2015,wu_neumann-michell_2017,garg_high-fidelity_2017,yang_integrated_2018,feng_multidisciplinary_2018,he_design_2019,miao_hull_2020,ni_multiple_2020,wang_many-objective_2021,villa_effective_2021,demo_efficient_2021,khan_regional_2021,demo_hull_2021,zhang_research_2021}, radial basis functions (RBF) \cite{yang_hydrodynamic_2014,yang_overview_2016,nazemian_multi-objective_2021,harries_application_2021,chang_dynamic_2021,zheng_dynamic_2021,nazemian_shape_2022}, arbitrary shape deformation \cite{zhang_research_2017,tezdogan_investigation_2018,zhang_hull_2018,nazemian_automated_2020,nazemian_cfd-based_2021}, patches \cite{tahara_cfd-based_2006,saha_hydrodynamic_2004,saha_hydrodynamic_2005,zhang_optimization_2018}, blending \cite{tahara_cfd-based_2006,hong_self-blending_2017,zong_hull_2018}, and morphing \cite{kandasamy_multi-fidelity_2011}, accounting for 28\%, are indicative of the need for more adaptable and flexible design approaches. Overall, Splines family (Spline, NURBS, B-Spline, T-Spline) approaches are the most used among the FPM, whereas FFD is the most used among the PPM methods. 

%(MAYBE ADD DETAIL FOR COMPERING SPLINE VERSUS FFD ...)

Figure \ref{fig:design_param}b illustrates the distribution of design space dimensionalities {and the cumulative sum of the associated occurrences}. Most studies concentrate on problems with 10 dimensions or fewer, indicating a focus on moderately complex design challenges. However, the presence of problems with higher dimensionality, greater than 50 \cite{percival_hydrodynamic_2001,klanac_optimization_2009,ehlers_procedure_2010,garg_high-fidelity_2015,garg_high-fidelity_2017,dagostino_design-space_2020}, up to 420 dimensions \cite{li_shape_2021}, reveals the presence of applications with highly complex and high-dimensional optimization challenges. {These high-dimensional optimizations are often facilitated by the use of adjoint gradients \cite{lee_surrogate_2016, garg_high-fidelity_2017, li_shape_2021, hamed_multi-objective_2022, nazemian_shape_2022}, since the computational cost of adjoint gradients scales favorably with the number of problem dimensions. Despite this success, adjoint solvers are not commonly used in the maritime field. This could be due to the relatively high complexity of these solvers which hampers a widespread adoption of the adjoint method for high-dimensional problems. Because of its high potential, research on adjoints for optimization should receive more attention.} {It is finally important to note that a significant portion of the works reviewed, approximately 26\%, do not explicitly specify the dimensionality of the design space. This omission indicates a gap in the reported information, meaning the presented distribution may only partially represent the problem dimensionalities encountered in SBDO research. The absence of detailed dimensionality data underscores a potential area for improvement in the clarity and completeness of reporting in the field.}

{The problem dimensionality} diversity raises the issue of the \textit{curse of dimensionality} \cite{Bellman:1957}, where larger design spaces exponentially increase computational costs and complicate the optimization process. Despite the variety of methods used for SBDO, considering both FPM and PPM, the definition of the design space still represents the true bottleneck in design processes. By limiting free variables, parametric models can significantly save time and costs. Hence, choosing restrictions based on experience, constraints from production, operational requirements, and market acceptance is crucial. Good parametric models stem from conscious choices of restriction, emphasizing the need for dimensionality reduction techniques in SBDO.

The development of dimensionality reduction techniques for shape optimization only recently gained attention. The simplest method to reduce the dimensionality of the design space is to identify the most important variables for the design problem and discard the remaining ones by setting them to a constant value during the optimization process, i.e. a factor screening, also known as feature selection. This process is conducted off-line (or upfront) the SBDO procedure. Sensitivity analysis has been used in \cite{zhang_sensitivity_2016} to prescribe the design space, whereas Pearson correlation coefficient has been used in \cite{mittendorf_hydrodynamic_2021} as a variable screening metric. On the contrary, online methods (during the SBDO procedure) have been proposed addressing dynamic space reduction in \cite{geremia_hull_2012,zheng_dynamic_2021}, where not the dimensionality of the design space is assessed, but the design variable range, exploring roughly the whole design space at the beginning of the SBDO and then restricting the variables range runtime, focusing on the most interesting part of the domain.
However, these approaches do not always provide the best solution, since factor screening is not able to evaluate the importance that the fixed variables could have during the optimization process, especially when combined with other variables, and dynamic space reduction could not take into account possible multi-modalities of the objective function, thus missing the optimum region. 
Hence, industrial design, in general, is increasingly searching for such dimensionality reduction methods that can capture, in a reduced-dimensionality space (possibly upfront), the underlying most promising directions of the original design space, preserving its relevant features and thereby enabling an efficient and effective optimization in the reduced space. The remedy has been found in dimensionality reduction techniques such as unsupervised learning, feature extraction, and modal representation, overall known as representation learning. These methods are capable of learning relevant hidden structures of the original design-space parameterization and have been developed focusing on the assessment of design-space variability and the subsequent dimensionality reduction before the optimization is performed.  
A method based on the Karhunen–Lo{\`e}ve expansion (KLE, equivalent to the proper orthogonal decomposition, POD) has been formulated in \cite{diez_design-space_2015} for the assessment of the shape modification variability and the definition of a reduced-dimensionality global model of the shape modification vector. No objective function evaluations nor gradients are required by the method. The KLE is applied to the continuous shape modification vector, requiring the solution of an eigenvalue problem for a Fredholm integral equation. The discretized Fredholm equation can be solved using principal component analysis. The method has been successfully applied to the optimization of the Delft catamaran in deterministic \cite{chen_high-fidelity_2015,serani_parameter_2016} and stochastic \cite{diez_stochastic_2018,pellegrini_hybridization_2020} conditions, the DTMB 5415 model \cite{serani_ship_2016}, Wigley hull \cite{liu_linear_2021}, as well as on different propellers \cite{gaggero_reduced_2020,doijode_machine_2022,doijode_machine_2022-1}.
Off-line methods improve shape optimization efficiency by reparameterization and dimensionality reduction, providing the assessment of the design space and the shape parameterization before optimization and/or performance analysis is carried out. The assessment is based on the geometric variability associated with the design space, making the method computationally very efficient and attractive (no simulations are required). Nevertheless, if the dimensionality reduction procedure is fed only with information on the shape modification vector, they may overlook the correlation between geometric variance and the actual objective function, since small variations in the geometry can produce significant variations in the objective function, e.g. flow separations and cavitation. For this reason, dimensionality reduction based on KLE has been extended to include physical information related to the optimization problem, resulting in significant improvements in both deterministic \cite{serani_adaptive_2019,pellegrini_derivative-free_2022} and stochastic \cite{serani_hull-form_2022} cases. A similar approach has been achieved via the active subspace method \cite{khan_regional_2021,demo_hull_2021}, which involves the identification of the so-called active subspaces of the input parameter space by analyzing the sensitivity of the output with respect to the input parameters, often using gradient information. Obviously, the use of physical information has a computational cost and cannot always be afforded by designers upfront the SBDO procedure. For this reason, a further attractive proposal is to substitute physical information with physics-related geometrical parameters. A recent example has been provided in \cite{khan_geometric_2022} where geometric moments are used to include physics information, applying it to two different ships.

%=================================================================================
\subsection{Numerical Solvers}
\begin{figure}[!t]
    \centering
    \includegraphics[width=0.5\columnwidth]{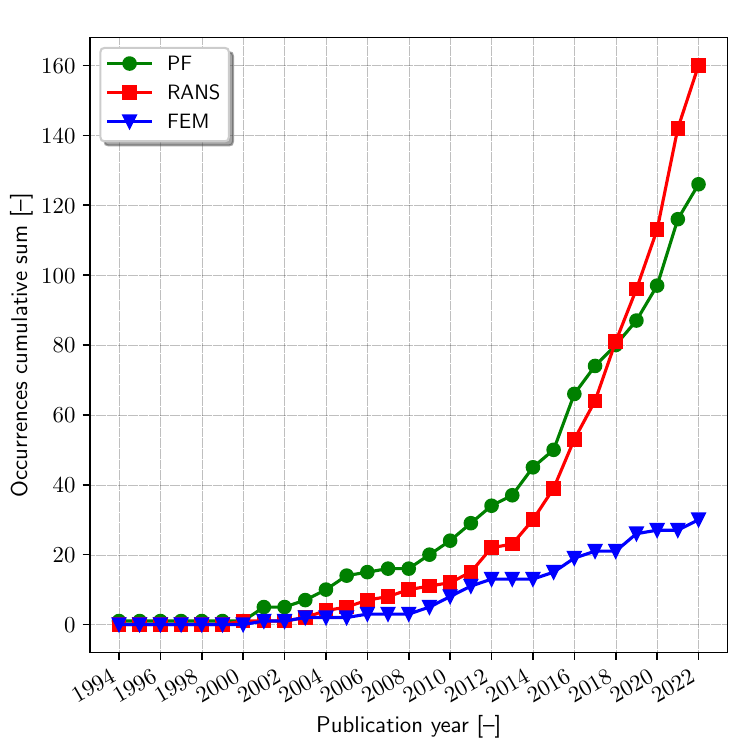}
    \caption{Cumulative sum of the kind of solvers used as a function of the publication year.}
    \label{fig:solvers}
\end{figure}
Figure \ref{fig:solvers} presents a compelling overview of the evolving solver usage in SBDO studies from 1994 to 2022. The graph shows the cumulative sum of occurrences for various solvers. These are potential flow methods (PF), Reynolds-averaged Navier-Stokes (RANS), and the finite element method (FEM). Each solver represents distinct computational approaches in SBDO.

The PF solver, while exhibiting a consistent increase in cumulative occurrences over the years, has been outpaced by the RANS solver since 2018. The increase in PF usage indicates its continued relevance, particularly in problems where potential flow assumptions are valid, such as in the early stages of aerodynamic or hydrodynamic design. PF solvers are mainly based on the boundary elements method (BEM), see e.g. \cite{campana_new_2009,brizzolara_concept_2012,danisman_reduction_2014,muratoglu_design_2017,li_design_2019,gaggero_reduced_2020,tao_optimized_2021}, but other examples have been found, such as strip theory \cite{grigoropoulos_hull_2004,sarioz_inverse_2009,bagheri_genetic_2014,bonfiglio_improving_2018,subramanian_genetic_2020}, slender body \cite{nazemian_global_2021}, vortex lattice \cite{vesting_surrogate_2014,vesting_development_2016}, and blade element momentum \cite{kamarlouei_multi-objective_2014,favacho_contribution_2016} methods, as well as isogeometric analysis combined with BEM \cite{kostas_ship-hull_2015,kostas_shape-optimization_2017,khan_geometric_2022}. It is important to recognize that within the realm of PF solvers, a significant portion are developed as proprietary, in-house tools, tailored to specific research or industrial needs. This trend underscores the specialized nature of PF solvers, which often require customization to address unique challenges in fluid dynamics and hydrodynamics. Nevertheless, commercially available options have also been used, see e.g. \cite{chen_inverse_2004,zhang_sensitivity_2016,chen_optimization_2016,cheng_hull_2018,liu_cross-entropy_2021,zheng_application_2021,harries_application_2021,chang_dynamic_2021,han_changwan_optimal_2014,bao_parametric_2021,nazemian_multi-objective_2022}.
%
%, SHIPFLOW (see e.g. \cite{chen_inverse_2004,zhang_sensitivity_2016,chen_optimization_2016,cheng_hull_2018,liu_cross-entropy_2021,zheng_application_2021,harries_application_2021,chang_dynamic_2021}) and ANSYS-AQWA \cite{han_changwan_optimal_2014,bao_parametric_2021,nazemian_multi-objective_2022} represent two of the few commercially developed PF solvers. 

The RANS solver shows a quartic trend in its cumulative occurrences. This significant rise reflects the growing preference for RANS in SBDO studies. The main cause is likely due to its enhanced capability in capturing complex turbulent flows and its applicability in a broader range of fluid dynamics problems compared to PF. This, in combination with an increase of computational resources which makes RANS affordable for practical applications, results likely in a strong increase of RANS usage over the years. The quartic nature of the trend suggests an accelerating adoption rate, highlighting RANS as an increasingly preferred tool for fluid dynamics optimization in recent years, as also reflected by the distribution between commercial (see, e.g., \cite{poloni_hybridization_2000,cirello_numerical_2008,mahmood_computational_2012,joung_shape_2012,vasudev_multi-objective_2014,chrismianto_parametric_2014,shi_optimal_2015,leifsson_optimal_2015,du_design_2016,gao_hull_2016,huang_optimization_2016,alam_design_2017,mizzi_design_2017,halder_wave_2018,zhang_computational_2018}), in-house developed \cite{duvigneau_hydrodynamic_2004,tahara_computational_2008,renaud_multi-objective_2022}, and open-source \cite{vernengo_physics-based_2016,berrini_geometric_2017,simonetti_optimization_2017,guerrero_surrogate-based_2018,coppede_hydrodynamic_2019,wang_investigation_2019,abdollahzadeh_optimization_2020,wang_three-dimensional_2022} solvers that is notably balanced. Commercial tools are widely used in various industries for their comprehensive capabilities and robust support structures. On the other hand, there are several notable in-house RANS solvers, which are developed within academic or research institutions for specific applications or research purposes. %Additionally, OpenFOAM stands out as a unique option, being both an in-house developed and open-source solver, widely used in the CFD community for its versatility and accessibility.

Finally, the use of FEM solvers \cite{he_design_2019} shows a more limited cumulative occurrence in SBDO studies despite its critical role in structural analysis. This might be indicative of the specific focus of the studies under consideration, possibly skewed more towards fluid dynamics than structural optimization. However, the presence of FEM, mainly composed of commercial software, see e.g., \cite{ehlers_particle_2012,sun_optimal_2012,dong_kriging-based_2017,jia_design_2019}, underscores its importance in the SBDO landscape, particularly for problems involving structural response and material optimization.

\begin{figure*}[!b]
    \centering
    \subfigure[Trend]{\includegraphics[width=1\textwidth]{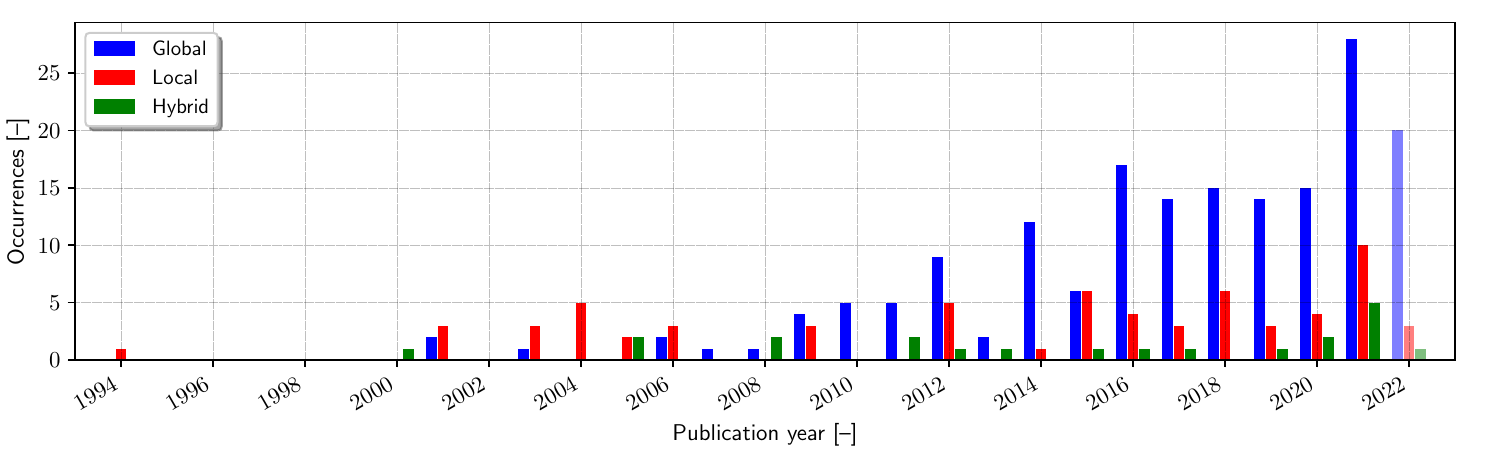}}\\
    \subfigure[Global]{\includegraphics[width=0.495\textwidth]{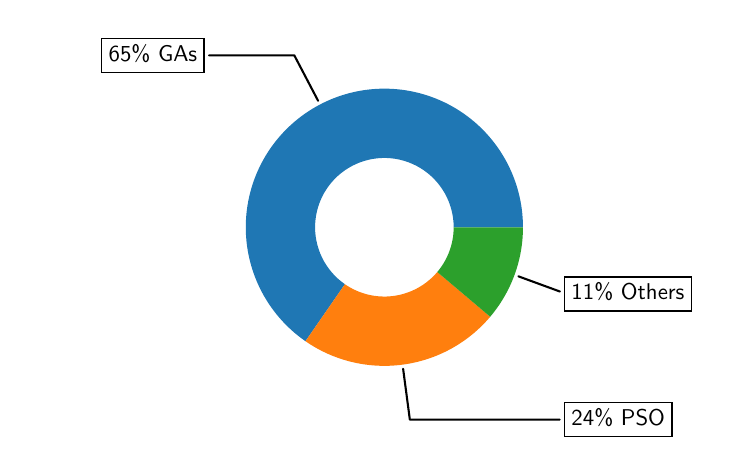}}
    \subfigure[Local]{\includegraphics[width=0.495\textwidth]{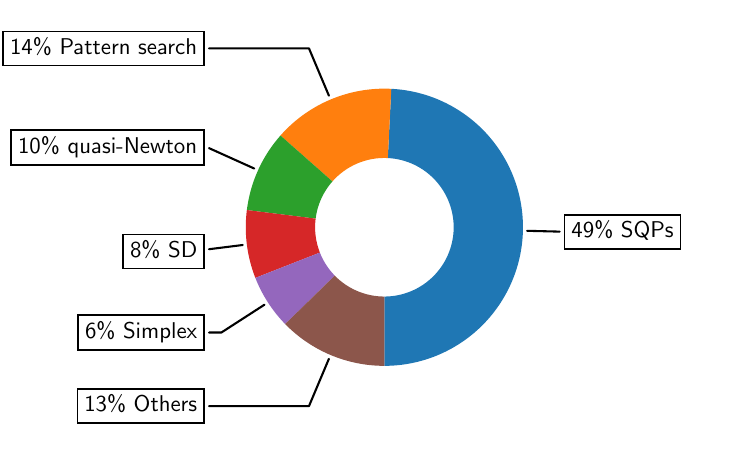}}
    \caption{Optimization algorithm occurrences (a) trend and subdivision by (b) global and (c) local categories.}
    \label{fig:opt_category}
\end{figure*}

The trends observed in Fig. \ref{fig:solvers} are indicative of the evolving preferences and technological advancements in the field of SBDO. The overtaking of PF by RANS in recent years points to a paradigm shift in solver selection, driven possibly by the increasing complexity of design problems and the need for more sophisticated fluid dynamics modeling capabilities. The limited but present use of FEM highlights the diverse range of optimization challenges addressed in SBDO, necessitating a variety of computational tools to cater to different aspects of marine engineering design.

%=================================================================================

\subsection{Optimization Methods}
In the evolving landscape of SBDO, the selection of optimization algorithms and the possible integration of surrogate methods play pivotal roles. These strategies are key in navigating the complex design spaces and computational challenges inherent in SBDO. The choice between global, local, or hybrid algorithms, as well as the adoption of surrogate-based approaches versus surrogate-free methods, reflects a strategic balance between exploration and exploitation, accuracy, and computational efficiency.

%=================================================================================
\subsubsection{Algorithms}

Figure \ref{fig:opt_category}a illustrates the year-by-year usage of global, local, and hybrid algorithms in SBDO studies. The trend towards global optimization algorithms signifies a strategic shift in SBDO. Global algorithms, known for their ability to explore the entire design space, are increasingly favored. This preference likely stems from their stochastic nature and heuristic methods, which are adept at avoiding local optima: a critical advantage in complex, multimodal design landscapes. The rising trend of global algorithms suggests an industry-wide acknowledgment of the complexity and unpredictability inherent in SBDO problems.

Within the realm of global optimization, genetic algorithms (GAs, see, e.g.,  \cite{dejhalla_application_2001,kitamura_optimization_2003,jang_adaptive_2009,lee_shape_2010,brizzolara_automatic_2011,whitfield_collaborative_2012,lu_hydrodynamic_2016,liu_hull_2017,rotteveel_inland_2017,esmailian_systematic_2017,wang_seakeeping_2018,fu_multi-objective_2018,vasudev_shape_2018,tahara_variable_2019,cheng_multi-objective_2019,vasudev_multi-objective_2019,luo_parametric_2019}) and particle swarm optimization (PSO, see, e.g., \cite{korgesaar_assessment_2010,de_pina_tailoring_2011,zhang_optimization_2016,yin_hydrodynamic_2021,li_application_2022}) dominate. As shown in Fig. \ref{fig:opt_category}b, GAs cover 65\% of global methods, leveraging mechanisms inspired by biological evolution, such as selection, crossover, and mutation. This allows for a robust exploration of the design space, making them particularly effective for non-linear, discrete, or mixed-variable optimization problems. PSO, with 24\%, employs a swarm intelligence approach that simulates social behavior patterns, providing a balance between exploration and exploitation in the search process. Within the remaining 11\% of the global methodologies, several notable algorithms have been identified and warrant mention. These include the infeasibility-driven evolutionary algorithm \cite{mohamad_ayob_uncovering_2011,alam_design_2014,alam_design_2017}, simulated annealing \cite{mohamad_ayob_uncovering_2011,yang_robust_2015}, artificial bee colony \cite{huang_hull_2016,yang_overview_2016}, and dividing rectangles \cite{campana_derivative-free_2016,serani_ship_2016}.

\begin{figure*}[!b]
    \centering
    \subfigure[Trend]{\includegraphics[width=1\textwidth]{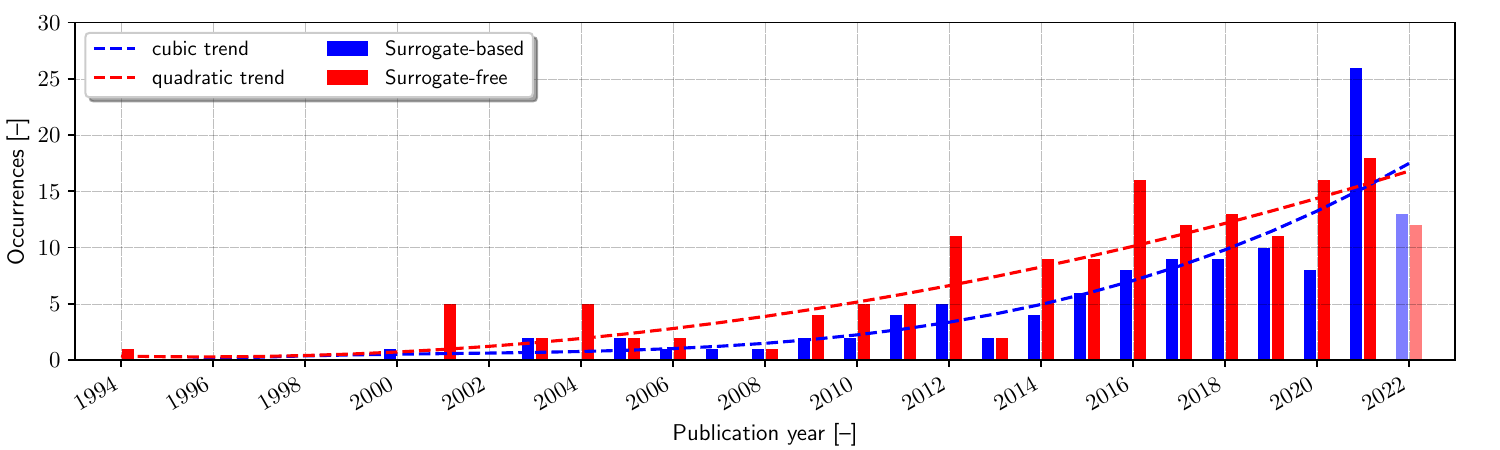}}\\
    \subfigure[Categories]{\includegraphics[width=0.495\textwidth]{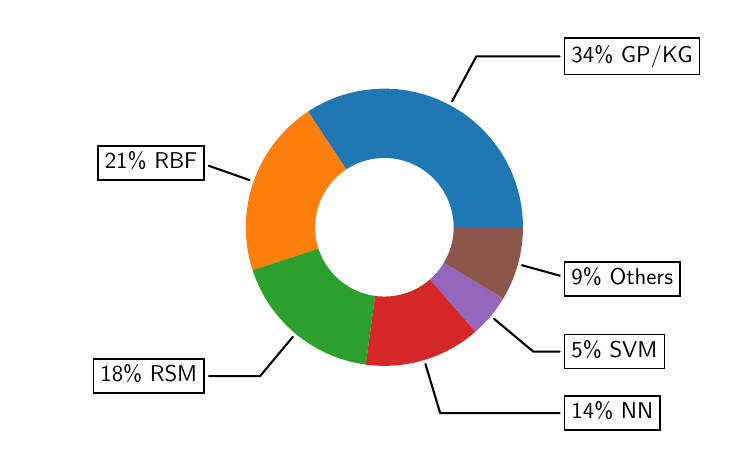}}
    \subfigure[DoE]{\includegraphics[width=0.495\textwidth]{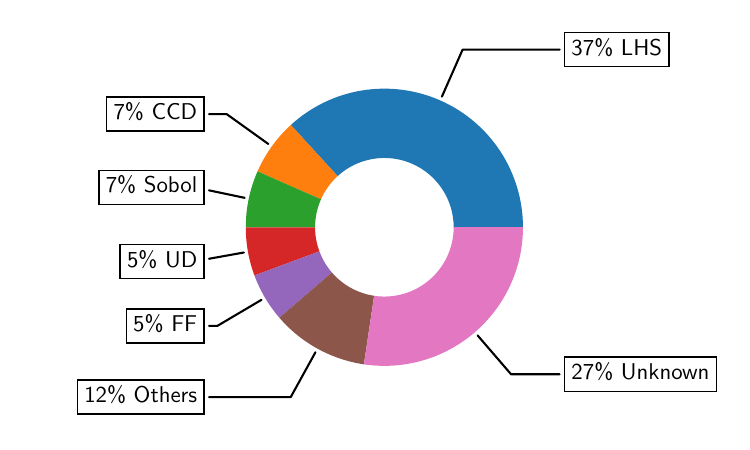}}
    \caption{Surrogate-based versus surrogate-free occurrence trends (a), surrogates categories (b), and design of experiments used for initial training (c).}
    \label{fig:surrogare_trend}
\end{figure*}
Considering local methods, the preference for sequential quadratic programming (SQP, see, e.g.,  \cite{cinquini_design_2001,tahara_computational_2004,suzuki_studies_2005,choi_hull-form_2015,baoji_research_2020,kunasekaran_design_2021,park_hull_2022}) and methods like quasi-Newton \cite{lowe_automatic_1994} methods (e.g., the Broyden-Fletcher-Goldfarb-Shanno, BFGS algorithm \cite{bonfiglio_improving_2018,bonfiglio_multi-fidelity_2020}) and pattern search, also known as Hooke and Jeeves algorithm \cite{grigoropoulos_hull_2004,sarioz_inverse_2009,sarioz_minimum_2012,du_design_2016}, as seen in Fig. \ref{fig:opt_category}c, aligns with problems where a good initial guess is available, and the design space is less rugged. In particular, SQP, with its ability to handle nonlinear constraints efficiently, is apt for fine-tuning solutions within a well-defined local region, complementing the global search methodologies. The steepest descent (SD) algorithm \cite{peri_design_2001}, the simplex method, also known as Nelder-Mead algorithm \cite{percival_hydrodynamic_2001,harries_application_2021,tran_optimal_2021,tran_optimization_2022}, and other gradient-based approaches \cite{kroger_adjoint_2018} are overall less preferred.

Finally, hybrid approaches deserve some hints. It may be noted that hybrid approaches include both memetic approaches (hybrid global/local) \cite{poloni_hybridization_2000,tahara_computational_2008,cirello_numerical_2008,vasile_ship_2013,yu_numerical_2015,pellegrini_hybridization_2020,luo_hull_2021,pellegrini_derivative-free_2022}, as well as hybridization of different global algorithms \cite{tahara_single-_2011}, global methods with reinforcement learning \cite{chiong_reinforcement_2012}, and local algorithms with multi-start approaches \cite{dong_kriging-based_2017,li_design_2019}. Among the memetic approaches the SHERPA (simultaneous hybrid exploration that is robust, progressive, and adaptive) algorithm \cite{ge_optimization_2016,im_duct_2020,nazemian_multi-objective_2021,nazemian_cfd-based_2021,nazemian_automated_2020,nazemian_simulation-based_2023}, noted for its robust and adaptive capabilities in handling complex design challenges, is gaining recognition in various engineering domains, not only marine. However, its proprietary nature, being exclusive to a specific software environment, presents potential limitations in terms of widespread adoption and accessibility, particularly in academic and open-source research communities where transparency and adaptability of algorithms are often paramount.

%=================================================================================
\subsubsection{Surrogates}

Figure \ref{fig:surrogare_trend}a compares the trend of solving SBDO problems with and without surrogate methods. The recent overtaking of surrogate-based methods over surrogate-free approaches marks a significant development in SBDO. 
In surrogate-based optimization, the original optimization problem in Eq. \ref{eq:soform} is reformulated by approximating the objective function ${f}(\bfx)$ and the eventual functional constraints ${g}_i(\bfx)$ and $h_j(\bfx)$ with surrogate models, denoted as $\hat{f}(\bfx)$, $\hat{g}_i(\bfx)$, and $\hat{h}_j(\bfx)$ respectively. This approach transforms the original optimization task into a more computationally tractable form by minimizing the surrogate objective function while ensuring that surrogate constraints are satisfied. The reformulated optimization problem is expressed as:
\begin{align}\label{eq:soform}
\min_{\bfx} \quad & \hat{f}(\bfx,\bfy) \\
\text{subject to} \quad & \hat{g}_i(\bfx,\bfy) \leq 0, \quad i = 1, \ldots, m \nonumber \\
& \hat{h}_j(\bfx,\bfy) = 0, \quad j = 1, \ldots, p . \nonumber
\end{align} 
Surrogate models, serving as approximations of the actual objective and constraint functions, offer substantial computational savings. The cubic trend of surrogate-based methods (see Fig. \ref{fig:surrogare_trend}a) reflects their growing importance in dealing with high-fidelity simulations that are computationally expensive, allowing for more iterations and a deeper exploration within feasible turnaround times.

The predominance of Gaussian process (GP, see, e.g.,  \cite{bonfiglio_improving_2018,coppede_hydrodynamic_2019,mittendorf_hydrodynamic_2021,demo_efficient_2021,renaud_multi-objective_2022,doijode_machine_2022}) and Kriging (KG, see, e.g., \cite{he_multidisciplinary_2011,tahara_single-_2011,xiao_generalised_2012,lee_metamodel-based_2012,yang_robust_2015,lee_surrogate_2016,liu_hull_2017,wu_neumann-michell_2017,nouri_optimization_2018,guerrero_surrogate-based_2018,huang_multi-objective_2019,thandayutham_optimization_2019,miao_cfd-based_2020,liu_hull_2021,cairns_numerical_2021,liu_optimization_2021,liu_resistance_2022}) methods (34\%) in surrogate-based optimization, as shown in Fig. \ref{fig:surrogare_trend}b, underscores their efficacy in capturing complex, nonlinear relationships with a relatively small number of samples. When it comes to practical applications, the distinction between GP models and KG models can become blurred despite their differences in original contexts and typical interpretations. This is particularly true in the context of surrogate modeling. In many cases, especially in computer experiments and design of experiments, the terms are used interchangeably, as the underlying mathematical principles are very similar. Both methods are highly appreciated for their ability to provide accurate predictions (excelling in modeling smooth functions) and a statistical framework that quantifies prediction uncertainty which is crucial for decision-making in the optimization process. However, computational challenges occur when applied to large datasets. Other popular methods like RBF (21\%, see, e.g., \cite{peri_high-fidelity_2005,tahara_computational_2008,chen_high-fidelity_2015,huang_hull_2016,zhang_sensitivity_2016,qiu_multi-objective_2019,wu_robust_2020,harries_application_2021,nazemian_automated_2020}), response surface methodologies (RSM, 18\%, see, e.g., \cite{duvigneau_role_2003,yang_metamodel_2011,li_multiobjective_2013,chrismianto_parametric_2014,luo_application_2015,shi_optimal_2015,huang_optimization_2016,hong_self-blending_2017,zong_hull_2018,lin_automatic_2019,yang_improved_2020,sun_design_2021,guan_parametric_2021-1}), neural networks (NN, 14\%, see, e.g., \cite{poloni_hybridization_2000,besnard_constructive_2007,hefazi_multidisciplinary_2010,grigoropoulos_hull-form_2010,vesting_surrogate_2014,han_changwan_optimal_2014,danisman_reduction_2014,sun_prediction_2019,ambarita_computational_2021,wang_three-dimensional_2022}), and support vector machines (SVM, 5\%, see, e.g., \cite{sun_optimal_2012,feng_multidisciplinary_2018,guo_cfd-based_2020,feng_multi-objective_2022,lv_optimisation_2022}) each offer unique advantages, such as local approximation capabilities and flexibility in modeling complex patterns. 
Specifically, RBFs are beneficial for multidimensional interpolation and smooth transitions, though they can struggle with larger, high-dimensional data; RSM is effective for design of experiments and process optimization but is less suited for non-linear or complex problems and requires extensive experimentation for accurate modeling; NNs, with their flexibility for complex relationships, are ideal for large datasets, but require significant data and are computationally intensive; lastly, SVM provide robust performance in high-dimensional spaces but are sensitive to kernel and parameter choices and computationally demanding for large datasets. This nuanced understanding of each method's strengths and weaknesses is crucial in guiding the selection of the most appropriate surrogate modeling technique for specific engineering optimization problems. Finally, the other 9\% surrogate-based approaches are composed of trust-region methods \cite{peri_multidisciplinary_2003,campana_shape_2006}, elliptic basis function \cite{guan_parametric_2021}, orthogonal polynomial methods \cite{yang_data-driven_2022}, and hyper-surrogate approaches \cite{halder_wave_2018}, where multiple surrogate methods are used, like RSM, RBF, and KG, and then averaged to get the objective prediction.  

It may be noted that in the present scoping review, and under the statistics provided in Fig. \ref{fig:surrogare_trend}b, works that characterize RBF models as single-layer NN are categorized under the use of RBF surrogates, rather than as conventional NN implementations. This classification stems from the mathematical alignment of single-layer RBF networks with RBF interpolation, highlighting their role as surrogate modeling techniques. In these instances, the RBF's function is used primarily to approximate complex nonlinear relationships within the data, distinguishing it from the multi-layered, deep-learning frameworks typically associated with NNs.

\begin{figure*}[!b]
    \centering
    \includegraphics[width=1\textwidth]{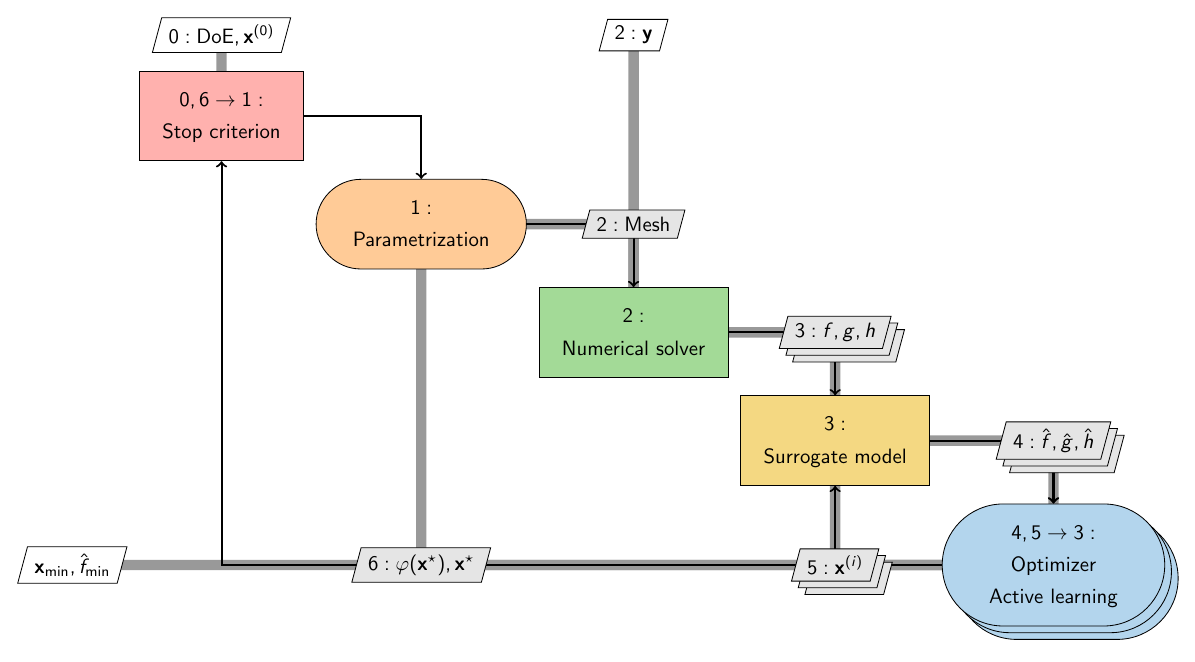}
    \caption{Example of extension of the XDSM diagram towards single-fidelity surrogate-based SBDO with active learning.}
    \label{fig:xdsm-surr}
\end{figure*}
Transitioning to another critical aspect of surrogate-based optimization, it is essential to acknowledge the pivotal role of the initial training and sampling approach employed for the surrogate models. The effectiveness of surrogate methods, as discussed earlier, hinges significantly on the quality and representativeness of the initial training data or design of experiments (DoE) used to construct these models. This data fundamentally influences the surrogate's ability to accurately capture the underlying behavior of the objective function and constraints. Therefore, the selection of an appropriate DoE becomes a key determinant in the success of surrogate-based optimization processes.
Among the various DoE employed (see Fig. \ref{fig:surrogare_trend}c), the Latin hypercube sampling (LHS), see, e.g., \cite{besnard_constructive_2007,yang_metamodel_2011,xiao_generalised_2012} covers 37\% of the cases (including optimal \cite{lee_metamodel-based_2012,yang_robust_2015} and universal \cite{coppede_hydrodynamic_2019} LHS) and this can be attributed to its effectiveness in generating well-distributed samples across the design space, ensuring a representative and unbiased training set for surrogate models. Other techniques include central composite design (CCD, 7\%, e.g. \cite{joung_shape_2012,chrismianto_parametric_2014,zong_hull_2018,serani_adaptive_2019,nazemian_simulation-based_2023}, Sobol (7\%, e.g. \cite{yang_improved_2020,villa_effective_2021,harries_application_2021}), uniform design (UD, 5\%, e.g. \cite{peri_multidisciplinary_2005,sun_optimal_2012,wang_aerodynamic_2021}), full factorial (FF, 5\%, e.g. \cite{thandayutham_optimization_2019,xu_optimization_2021,kunasekaran_design_2021}), and finally the remaining 12\% includes orthogonal arrays (5\%, \cite{peri_high-fidelity_2005,tahara_computational_2008,campana_numerical_2009}), Hammersley/Halton sequences \cite{chen_high-fidelity_2015,leotardi_variable-accuracy_2016,doijode_machine_2022}, as well as random/Monte Carlo sampling \cite{wang_investigation_2019,du_identification_2022,wang_three-dimensional_2022}. However, it is noteworthy that in 27\% of the cases, the specific DoE strategy employed remains unidentified or unspecified. This lack of clarity on the training approach used can have implications for the interpretability and reproducibility of the optimization results. Consequently, this highlights a gap in the current body of research, underscoring the need for more transparent and detailed reporting of the sampling methodologies in surrogate-based optimization studies to better understand their impact on the effectiveness of the surrogate models.

In the domain of surrogate-based optimization, the development of multi-fidelity or variable-fidelity methods has emerged as a key strategy to enhance the effectiveness of surrogate models while also conserving computational resources. These methods leverage varying levels of model fidelity, combining computationally expensive high-fidelity simulations with less costly lower-fidelity approximations in order to construct more informed and efficient surrogates. Despite their potential benefits, the scoping review reveals that only 12\% \cite{peri_multidisciplinary_2003,campana_shape_2006,tahara_computational_2008,campana_numerical_2009,kandasamy_multi-fidelity_2011,koziel_simulation-driven_2012,jiang_modified_2016,bonfiglio_improving_2018,serani_adaptive_2019,thandayutham_hydrostructural_2020} \cite{bonfiglio_multi-fidelity_2020,gaggero_marine_2022,liu_multi-fidelity_2022} of surrogate-based approaches have employed multi-fidelity methodologies, and their application appears sporadic over the years covered by the review. This limited utilization raises questions about the popularity and perceived benefits of multi-fidelity methods in this specific field. It is unclear whether this lack of widespread adoption is due to a general underutilization of these methods in the industry, or if there exist ambiguities and uncertainties regarding the actual advantages of integrating multi-fidelity approaches in surrogate-based optimization for marine engineering applications. This observation points to a potential area for further investigation and clarification, as the effective use of multi-fidelity methods could significantly impact the efficiency and accuracy of optimization processes in this domain.

In concluding the discussion on surrogate-based optimization, it is crucial to recognize the role of adaptive sampling or active learning methods in enhancing the effectiveness of these models. Such techniques, for both single- and multi-fidelity methods, start with an initial DoE, subsequently adapted by incorporating new samples $\bfx^\star$ in areas most beneficial for optimization. A variety of strategies have been employed for this purpose, including, among others, the so-called acquisition function $\varphi$ based on: the validation of the best found \cite{peri_multidisciplinary_2005,grigoropoulos_hull-form_2010,tahara_single-_2011,koziel_simulation-driven_2012,chrismianto_parametric_2014,chen_high-fidelity_2015,serani_ship_2016,guerrero_surrogate-based_2018,zhang_optimum_2019,cairns_numerical_2021}, the maximum uncertainty \cite{campana_numerical_2009,serani_adaptive_2019}, the expected improvement \cite{lee_surrogate_2016,bonfiglio_improving_2018,serani_adaptive_2019,bonfiglio_multi-fidelity_2020}, and lower confidence bounding \cite{serani_adaptive_2019,luo_hull_2021}. These methods aim to iteratively refine the surrogate model by focusing on regions of the design space where additional information can significantly influence the optimization outcome. 
Despite the apparent advantages of these adaptive techniques, this scoping review indicates that in 21\% of the surrogate-based methods employing adaptive sampling approaches, the specific technique utilized remains unspecified. This lack of detail not only hinders the full understanding of the method's implementation but also obscures the comparative analysis of different techniques' efficacy. Given the potential impact of adaptive sampling on the accuracy and efficiency of surrogate-based optimization, particularly in marine engineering applications, this represents a significant gap in the current literature. A more transparent and detailed reporting of adaptive sampling methods could provide deeper insights into their benefits and limitations, fostering their more informed and effective use in the field. 

An example of how SBDO workflow shown in Fig. \ref{fig:xdsm} can be extended to the use of a general single-fidelity surrogate approach, including active learning, is given in Fig. \ref{fig:xdsm-surr}. 
{The diagram illustrates how the surrogate model acts as an intermediary between the numerical solver and the optimization algorithm. This arrangement facilitates the application of the optimization algorithm directly on the surrogate model to identify the optimal solution, denoted as 
$\bfx_{\min}$ and $\widehat{f}_{\min}$. Concurrently, an active learning-driven optimization procedure operates in parallel. This procedure employs an acquisition function, $\varphi$, to systematically pinpoint potential new candidate solutions $\bfx^{\star}$ to be sampled. These candidates are then processed through the numerical solver if the predefined stopping criterion has not yet been met. This dual-path approach integrates surrogate modeling with active learning to efficiently converge towards the optimum by balancing the exploration of the solution space and the exploitation of known high-potential areas.}
A further example of XDSM diagram extended to multi-fidelity methods can be found in \cite{spinosa2023simulation}.

%=================================================================================
%
\begin{figure}[!b]
    \centering
    \includegraphics[width=0.5\columnwidth]{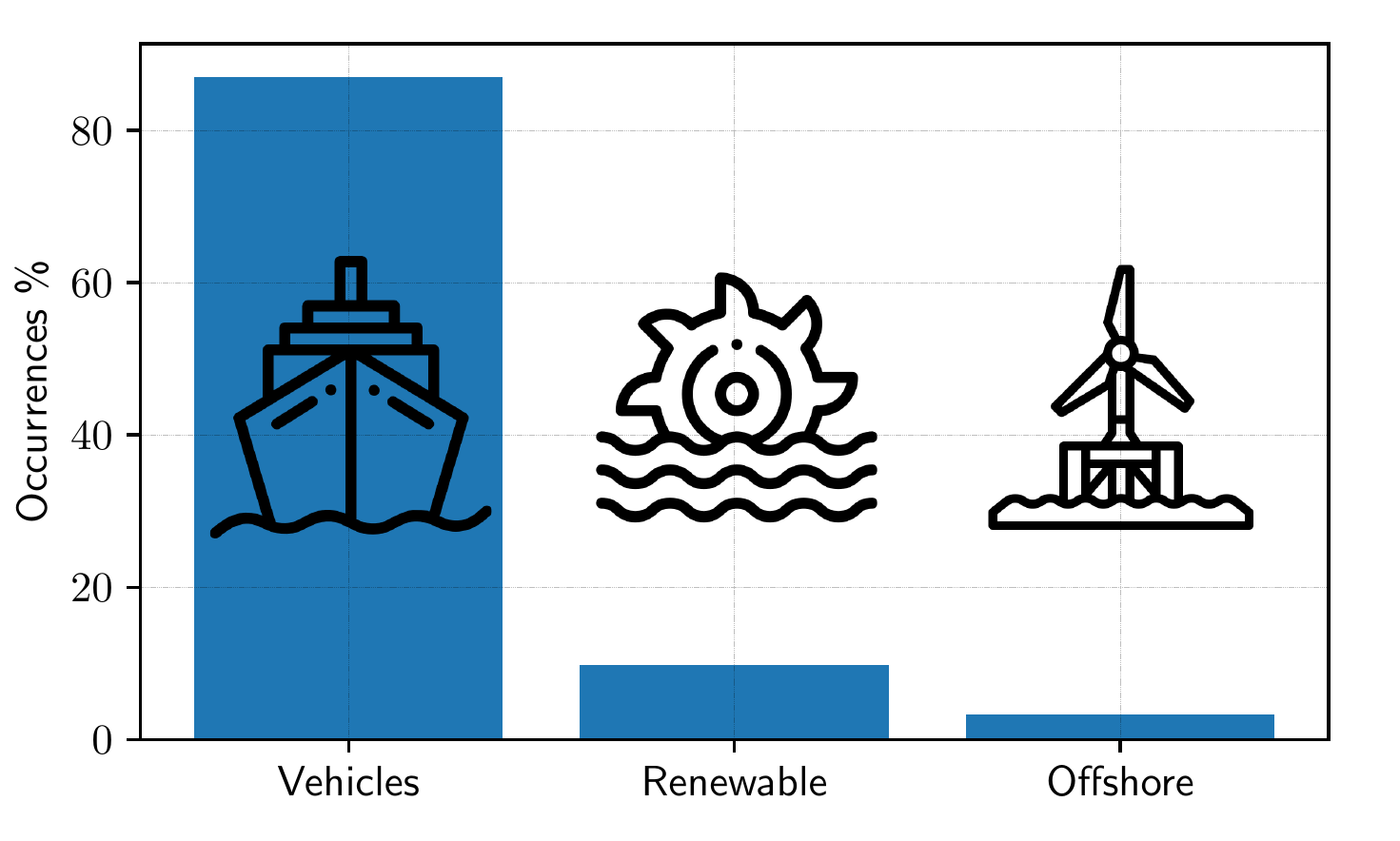}
    \caption{Occurrences of SBDO applied to marine engineering main applications fields.}
    \label{fig:chart-vehicles}
\end{figure}
\subsection{Applications}
Figure \ref{fig:chart-vehicles} shows the breakdown of the SBDO applications in marine engineering. The overwhelming majority of SBDO applications are dedicated to vehicle design (87\%), which includes ships (see, e.g., \cite{valorani_sensitivity_2003,zhang_optimization_2009,wilson_hull_2010,zhang_shape_2012,lv_trim_2013,dambrine_theoretical_2016}), submarines (see, e.g., \cite{ignacio_optimized_2019,page_simulation-driven_2020}), and various types of watercraft. This dominant focus can be attributed to several factors: (i) marine vehicles often have complex design requirements balancing hydrodynamic efficiency, stability, load capacity, and speed, consequently SBDO provides a powerful tool to optimize these competing factors; (ii) the marine vehicle industry is highly competitive, with a constant demand for improved performance and efficiency and SBDO enables designers to explore innovative shapes and configurations that might not be feasible through traditional design methods; (iii) the increasing environmental regulations and the push for energy efficiency drive the need for advanced optimization techniques to meet these stringent standards.
The use of SBDO in the development of renewable energy solutions in marine settings, such as wave \cite{coiro_diffuser_2016,e_silva_hydrodynamic_2016,simonetti_optimization_2017,halder_wave_2018,tao_optimized_2021,bao_parametric_2021} and ocean-thermal \cite{chen_optimal_2022} energy converters, pumps \cite{wang_multi-condition_2020}, and tidal \cite{kinnas_computational_2012,huang_optimization_2016,zhang_optimization_2016,sun_prediction_2019,im_duct_2020,khanjanpour_optimization_2020,ambarita_computational_2021,yeo_tidal_2022}, marine/ocean current \cite{yang_hydrofoil_2012,luo_multi-point_2014,thandayutham_optimization_2019,thandayutham_hydrostructural_2020,kunasekaran_design_2021,zhu_optimization_2022}, river hydrokinetic \cite{muratoglu_design_2017,barbaric_investigation_2020}, and offshore wind \cite{han_changwan_optimal_2014,du_design_2016,lemmer_semi-submersible_2020} turbines, highlights its growing importance, covering 10\% of the literature. This category's smaller proportion might be due to the relatively newer field compared to traditional marine vehicle design. Furthermore, the design of renewable energy systems involves complex interactions with the marine environment, requiring sophisticated models that can be challenging to optimize.
The smallest category in the breakdown is offshore applications (3\%), which include steel catenary risers \cite{yang_metamodel_2011,de_pina_tailoring_2011}, deep-sea test miners \cite{lee_metamodel-based_2012}, platforms and semi-submersible structures \cite{yang_robust_2015,qiu_multi-objective_2019,jang_fea_2019}, mooring systems \cite{li_design_2019}, and ocean bottom flying nodes \cite{xu_optimization_2021}. Factors influencing this lower percentage include high stakes and safety concerns, as well as complex environmental conditions. Offshore structures are often subject to stringent safety standards due to the high risks involved, possibly leading to a more conservative approach in adopting new optimization techniques. Moreover, the design of offshore structures must account for a wide range of environmental conditions, making the optimization process more challenging.

\begin{figure}[!t]
    \centering
    \includegraphics[width=0.5\columnwidth]{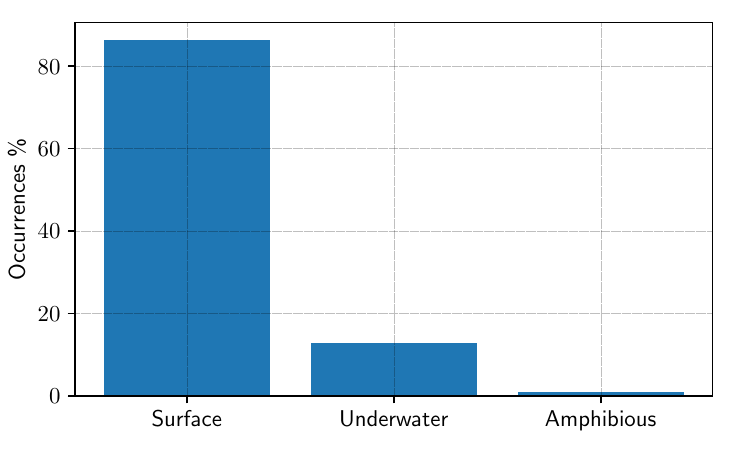}
    \caption{Occurrences of SBDO applied vehicles sub-categories.}
    \label{fig:marine-vehicles}
\end{figure}
Among vehicle design, Fig. \ref{fig:marine-vehicles} offers insights into where optimization efforts are being primarily focused. Specifically, 86\% is composed of surface vehicles, 13\% underwater, and the remaining 1\% amphibious.
%Surface Vehicles likely include a wide range of vessels such as cargo ships, passenger ships, ferries, and yachts. 
The optimization of surface vehicles can be pivotal in enhancing various aspects like hydrodynamic efficiency and seakeeping, resulting in less fuel consumption, improved stability and payload capacity. SBDO's significant role in surface vehicle design may be due to the large economic and environmental impact of these vessels, driving a need for continuous improvement in their performance and efficiency.
Underwater vehicles include submarines \cite{chrismianto_development_2015,vasudev_modular_2016} and autonomous underwater vehicles (AUVs, see, e.g., \cite{joung_shape_2012,vasudev_multi-objective_2014,alam_design_2014,gao_hull_2016}). The design optimization of these vehicles focuses on aspects like efficient maneuverability, stability under water, and energy efficiency for extended mission ranges. The application of SBDO in underwater vehicle design indicates a focus on specialized performance characteristics unique to the underwater environment, such as pressure resistance and stealth capabilities.
Finally, amphibious vehicles \cite{seth_amphibious_2020,du_identification_2022} are specialized vehicles that operate both in water and in air or land. The design challenges for amphibious vehicles are particularly complex due to the need to optimize performance in two very different environments. SBDO can play a key role in balancing these dual requirements, optimizing aspects such as buoyancy, stability, and propulsion efficiency.

\begin{figure}[!b]    \centering
    \includegraphics[width=0.5\columnwidth]{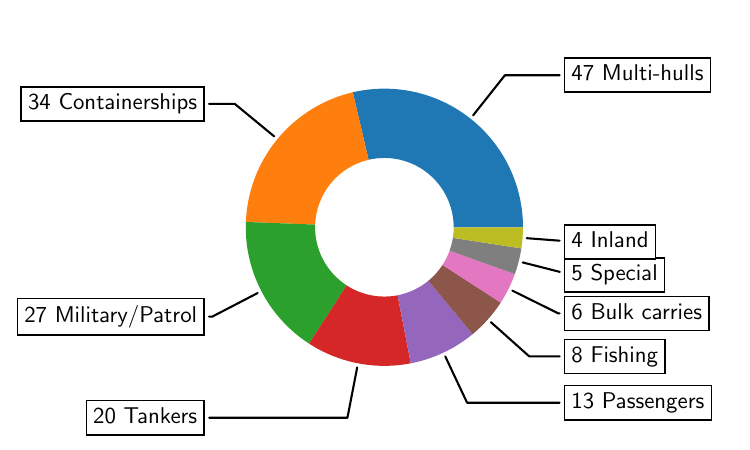}
    \caption{Occurrences of SBDO applied to surface vehicles sub-categories.}
    \label{fig:surface-vehicles}
\end{figure}
Due to the predominance of surface vehicles, a further breakdown has been conducted in this subfield. The sub-categories are shown in Fig. \ref{fig:surface-vehicles}.
A significant focus on containerships (see, e.g., \cite{kroger_adjoint_2018,feng_parametric_2021}) in SBDO applications aligns with their vital role in global trade. Optimization for these vessels likely focuses on maximizing cargo capacity, fuel efficiency, and minimizing environmental impact, crucial for cost-effective and sustainable operations. The Korea research institute of ships and ocean engineering (KRISO) container ship (KCS) represents the most used benchmark in this sub-category, see, e.g., \cite{vasile_ship_2013,chrismianto_parametric_2014,guha_amitava_application_2015,kroger_adjoint_2018,chen_parametric_2018,lu_hydrodynamic_2019,coppede_hydrodynamic_2019,miao_hull_2020,chang_dynamic_2021,khan_regional_2021,park_hull_2022}, serving as a standard reference model for various hydrodynamic studies.
The optimization of military \cite{peri_multidisciplinary_2005,zhang_minimum_2020,zha_hull_2021} and patrol \cite{mohamad_ayob_uncovering_2011,sarioz_minimum_2012} vessels underscores the importance of performance, stealth, and agility in these applications. SBDO can be instrumental in enhancing these attributes, contributing to the effectiveness and safety of naval operations. As for containerships, also military vessels have their specific standard benchmark, represented by the David Taylor model basin (DTMB) 5415 model, which has been extensively used for hull-form optimization purposes \cite{peri_multidisciplinary_2003,tahara_computational_2004,campana_shape_2006,tahara_computational_2008,li_multiobjective_2013,campana_derivative-free_2016,serani_ship_2016,wu_neumann-michell_2017,zhang_research_2017,zhang_optimization_2018,dagostino_design-space_2020,zha_hull_2021,serani_hull-form_2022,pellegrini_derivative-free_2022,liu_multi-fidelity_2022}.
The application of SBDO in tanker design (see, e.g., \cite{goren_mathematical_2017}) reflects the need for optimizing fuel efficiency and safety, given their role in transporting large volumes of liquid cargo, including oil and chemicals. The KRISO very large crude carrier (KVLCC2) model is the actual benchmark in this sub-category, see, e.g., \cite{duvigneau_role_2003,liu_optimization_2021}.
The application of SBDO in several further categories indicates a broad spectrum of optimization goals, from enhancing the efficiency of bulk carriers \cite{li_bow_2014,tahara_variable_2019,he_design_2019,liu_resistance_2022} and fishing \cite{leifsson_optimal_2015,tomasz_abramowski_energy_2017,hong_self-blending_2017,lin_hull_2018,yang_integrated_2018,tezdogan_investigation_2018,tran_optimal_2021,zhao_optimisation_2021} vessels to improving passenger comfort and safety in passenger's vessels \cite{chiong_reinforcement_2012,geremia_hull_2012,luo_design_2017,harries_application_2021}, including yachts \cite{lowe_automatic_1994,poloni_hybridization_2000,cirello_numerical_2008,diez_hydroelastic_2012,leotardi_variable-accuracy_2016,berrini_geometric_2017,bacciaglia_controllable_2021} and cruise ships \cite{wang_aerodynamic_2021,demo_efficient_2021}. The optimization of inland \cite{cinquini_design_2001,favacho_contribution_2016,rotteveel_inland_2017,maia_computational_2021} and special ships also points to specialized requirements, perhaps related to shallow waters navigation or unique operational roles like research vessels \cite{subramanian_genetic_2020,timurlek_hydrodynamic_2022}, survey ships \cite{li_application_2022}, or offshore aquaculture \cite{feng_multidisciplinary_2018,wang_many-objective_2021}.

\begin{figure}[!t]
    \centering
    \includegraphics[width=0.5\columnwidth]{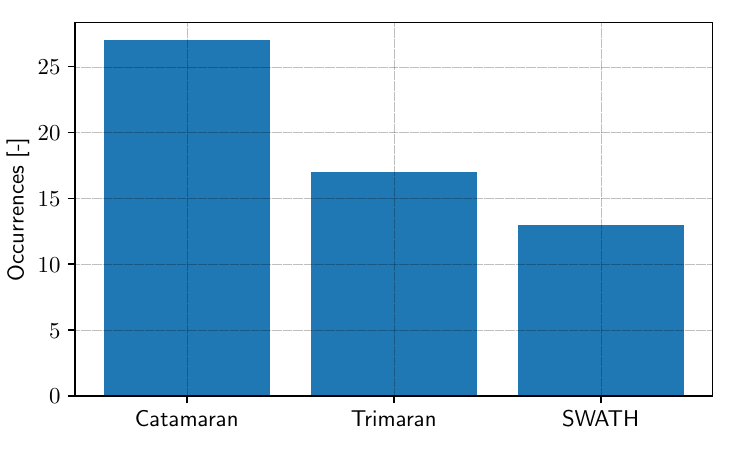}
    \caption{Occurrences of SBDO applied to multi-hulls sub-categories.}
    \label{fig:multi-hull}
\end{figure}
As shown in Fig. \ref{fig:surface-vehicles}, the strongest emphasis on surface vessels is represented by multi-hull designs, such as catamarans and trimarans, suggesting a focus on seakeeping and efficiency, resulting in improved stability and speed. Multi-hulls present unique design challenges that SBDO can help address, particularly in balancing stability with performance. For these reasons a deeper analysis has been conducted on multi-hull vessels, revealing three main sub-categories, which are catamarans, trimarans, and SWATH vehicles (see Fig. \ref{fig:multi-hull}).
Catamarans, with two parallel hulls of equal size, offer stability and spaciousness, making them popular for passenger ferries and recreational vessels. SBDO in catamaran design \cite{chen_inverse_2004,saha_hydrodynamic_2004} likely focuses on optimizing hull shape for stability \cite{campana_numerical_2009,tahara_single-_2011,diez_stochastic_2018} and reducing resistance, improving fuel efficiency \cite{kandasamy_multi-fidelity_2011,danisman_reduction_2014,miao_cfd-based_2020,mittendorf_hydrodynamic_2021,harries_application_2021,ahmad_fitriadhy_optimization_2022}. The standard benchmark model used for developing and assessing SBDO methodologies is represented by the Delft catamaran, see, e.g., \cite{diez_design-space_2015,chen_high-fidelity_2015,serani_parameter_2016,pellegrini_formulation_2017}. %The optimization might also focus on the aerodynamics of the vessel, crucial for high-speed applications.
Trimarans, featuring a main hull with two smaller outrigger hulls, are known for their speed and stability, making them suitable for high-speed ferries and racing yachts. In trimaran design \cite{hefazi_multidisciplinary_2010,yang_hydrodynamic_2014}, SBDO can play a crucial role in optimizing the hull configuration for balance and speed \cite{wang_seakeeping_2018,zong_hull_2018,nazemian_multi-objective_2022,nazemian_simulation-based_2023}, ensuring structural integrity \cite{jia_design_2019} while maximizing performance \cite{guo_cfd-based_2020,nazemian_multi-objective_2021,nazemian_automated_2020,nazemian_cfd-based_2021,nazemian_global_2021,wang_study_2021,hamed_multi-objective_2022,nazemian_multi-objective_2022-1,nazemian_shape_2022}.
The use of SBDO in trimarans can also address specific challenges like wave-piercing capabilities \cite{lv_optimisation_2022} and maneuverability, enhancing their performance in various marine conditions.
SWATH vessels are designed to minimize hull volume at the water's surface, reducing the impact of waves and providing a smoother ride in rough seas. SBDO in SWATH design \cite{besnard_constructive_2007,brizzolara_automatic_2011,xiao_generalised_2012,brizzolara_concept_2012,ni_multiple_2020,yang_improved_2020,bonfiglio_multi-fidelity_2020,guan_parametric_2021-1} is likely centered on optimizing the hull shape and configuration \cite{jiang_modified_2016,lin_automatic_2019,pellegrini_hybridization_2020} to achieve the desired stability and seakeeping qualities \cite{bonfiglio_improving_2018,renaud_multi-objective_2022}, making them ideal for applications like research vessels and coast guard ships.
It should be finally highlighted that Fig. \ref{fig:surface-vehicles} does not account for the hull-form studies applied to the Wigley \cite{percival_hydrodynamic_2001,dambrine_theoretical_2016} and systematic series S60 \cite{zhang_optimization_2009,mahmood_computational_2012,dejhalla_application_2001,dejhalla_numerical_2001,huang_hull_2016} benchmark models because they cannot be included in any of the specified subcategories. Nevertheless, they have been used for specific development/assessment of SBDO method \cite{zakerdoost_evolutionary_2013,bagheri_genetic_2014,bagheri_optimizing_2014,hou_ship_2017,zhang_computational_2018,zhang_hull_2018,zhang_minimum_2020,liu_linear_2021,zhang_research_2021,yang_overview_2016,zheng_application_2021,zheng_dynamic_2021,cheng_hull_2018,yin_hydrodynamic_2021}, as well as for particular operational/environmental conditions, like high speed \cite{baoji_research_2020} and shallow waters \cite{saha_hydrodynamic_2004}, or retrofitting \cite{liu_hull_2021,esmailian_systematic_2017}.

\begin{figure}[!b]
    \centering
    \includegraphics[width=0.5\columnwidth]{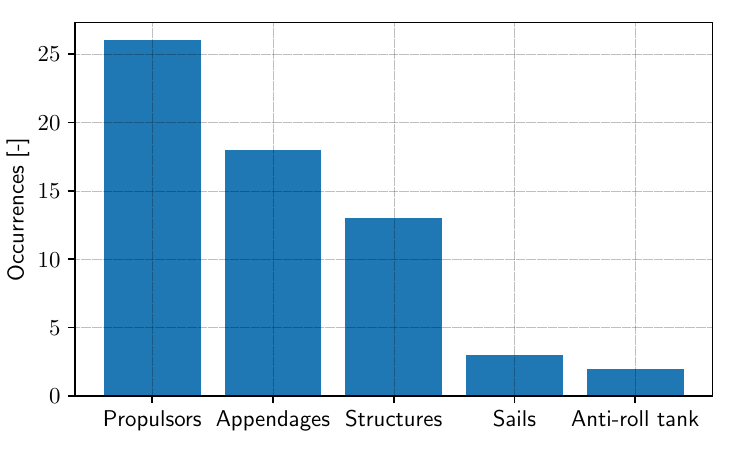}
    \caption{Occurrences of SBDO applied to marine components following the primary classification level.}
    \label{fig:app-detail}
\end{figure}
Finally, a breakdown of SBDO applied to marine components is presented in Fig. \ref{fig:app-detail}.
Propulsors, including propellers \cite{young_reliability-based_2010}, water jets \cite{guo_cfd-based_2020,huang_multi-objective_2019}, and thrusters \cite{feng_multi-objective_2022}, are critical for the movement and maneuverability of marine vehicles. Shape optimization in this area focuses on improving hydrodynamic efficiency \cite{whitfield_collaborative_2012,kamarlouei_multi-objective_2014,ma_design_2014,mizzi_design_2017,nouri_optimization_2018,gaggero_marine_2022,doijode_machine_2022}, reducing cavitation \cite{bertetta_cpp_2012,yu_numerical_2015,mirjalili_confidence-based_2018,gaggero_reduced_2020}, and minimizing noise \cite{vesting_surrogate_2014,vesting_development_2016,jiang_prediction_2020}. The optimization could involve refining blade shapes and angles \cite{lu_research_2021} to enhance propulsion efficiency while reducing fuel consumption \cite{bacciaglia_controllable_2021,doijode_machine_2022-1} and environmental impact \cite{favacho_contribution_2016,gaggero_design_2017,esmailian_systematic_2017}, including also retrofitting solutions, like equalizing ducts \cite{furcas_design_2020}.
Marine vehicle appendages include rudders \cite{chen_parametric_2018}, fins \cite{garg_high-fidelity_2015,garg_high-fidelity_2017}, and keels \cite{poloni_hybridization_2000}, which play essential roles in stability and steering. Shape optimization in appendages \cite{duvigneau_hydrodynamic_2004,vernengo_physics-based_2016,kostas_shape-optimization_2017,serani_adaptive_2019} aims to enhance hydrodynamic performance, improve maneuverability, and reduce drag \cite{wang_investigation_2019}. This might involve optimizing the size, shape, and positioning of these components to achieve a balance between stability and agility \cite{diez_hydroelastic_2012,leotardi_variable-accuracy_2016}.
Structures likely encompass the hull and superstructure of marine vehicles, as well as substructures of offshore platforms \cite{han_changwan_optimal_2014,yang_robust_2015,qiu_multi-objective_2019,jang_fea_2019,zhu_optimization_2022}. Shape optimization in structures focuses on enhancing overall hydrodynamic performance, maximizing space utilization, and ensuring structural integrity \cite{ha_design_2006,choi_sampling-based_2015,dong_kriging-based_2017,jia_design_2019}. In addition, it involves tweaking hull forms for better wave resistance, stability, and seakeeping qualities, crucial for efficiency and safety \cite{kitamura_optimization_2003,diez_experimental_2022}, including crashworthiness \cite{klanac_optimization_2009,ehlers_procedure_2010,ehlers_particle_2012,korgesaar_assessment_2010}.
In sailboats and sailing yachts, the optimization of sail shapes is vital for maximizing wind propulsion efficiency \cite{cairns_numerical_2021}. This involves determining the optimal curvature, material stiffness, and positioning of sails to harness wind power effectively, which is essential for performance in competitive sailing and leisure cruising \cite{lee_surrogate_2016,ma_hard_2019}.
Finally, anti-roll tanks are used to stabilize ships by reducing rolling motion caused by waves. Shape optimization in anti-roll tanks aims to maximize their effectiveness in damping roll motion while minimizing the impact on the vessel's overall performance and weight distribution \cite{subramanian_genetic_2020,liu_optimization_2021}.

The detailed breakdown of shape optimization in various marine vehicle components underscores the comprehensive and multifaceted nature of design challenges in marine engineering. Shape optimization in each of these areas requires a deep understanding of fluid dynamics, material properties, and operational conditions. The focus on specific components like propulsors, appendages, and structures reflects the industry's commitment to enhancing performance, safety, and environmental sustainability. The optimization of sails and anti-roll tanks highlights specialized areas where SBDO can significantly impact vessel performance and passenger comfort.
This analysis demonstrates the critical role of shape optimization in advancing the design and functionality of marine vehicles. It highlights the technological advancements in SBDO and its application in addressing the nuanced and complex design requirements of different components of marine vehicles.

Overall, this analysis underscores the adaptability and potential of SBDO across various facets of marine engineering, promising continued innovation and improvement in the design of marine vehicles, renewable energy systems, and offshore structures.

%=================================================================================

\section{Discussion}
The marine engineering field, while recognizing the advantages of more comprehensive multi-objective and stochastic optimization approaches, shows a marked preference for simpler, deterministic single-objective formulations. This trend results from the tendency to provide a simple and clear demonstration of new SBDO methodologies. At the same time, it highlights important areas for future growth such as the adoption of stochastic problem formulations, such as RDO, RBDO, and RBRDO. These approaches more accurately reflect the uncertainties characteristic of marine environments and align with broader trends in marine engineering, including digitalization, sustainability, and evolving regulatory landscapes. 
The analysis of problem formulations in SBDO studies reveals a landscape where constrained problems dominate, reflecting the complex nature of marine engineering challenges. However, the significant proportion of studies with unclear formulations and the apparent gap in the discussion of constraint-handling strategies highlight areas for improvement in SBDO research. Future studies would benefit from a more explicit focus on the nature and management of constraints, thereby enhancing the relevance, applicability, and impact of SBDO in marine engineering.
The scarcity of MDO applications also highlights a potentially huge area for growth and development in marine engineering research. As the field continues to develop, an increased recognition of the benefits of a more integrated multidisciplinary approach is expected. MDO is especially useful in tackling complex design challenges that encompass multiple engineering facets. Future research could focus on developing more accessible and efficient MDO methodologies, facilitating their broader adoption in marine engineering optimization problems.

The variety of parameterization techniques reflects a range of approaches to defining design spaces, while the distribution of design space dimensionalities reveals both a focus on more manageable problems and an interest in tackling more complex, high-dimensional optimization challenges. This analysis underscores the need for continued innovation in SBDO methodologies, particularly in addressing the challenges posed by high-dimensional design spaces, and overcoming the curse of dimensionality. Dimension reduction techniques such as factor screening, sensitivity analysis, and dynamic space reduction are classical approaches to mitigate the curse of dimensionality. However, these techniques do not capture multi-modalities of the objective function and may therefore fail to find the optimum region. Unsupervised learning, feature extraction, and representation learning such as KLE and POD overcome these issues and do not require objective function evaluations or gradients. These methods are based on geometrical variance and do not account for the relation between geometrical variation and the variation of the objective. The inclusion of physical (objective) information is therefore identified as a promising way to improve dimension reduction techniques. Nevertheless, for practical application in an industrial context, where parametrization methods are mainly CAD-based, designers cannot easily retrieve the original design variables from the reduced design space (also known as the pre-image problem). It can be noted that a back-mapping procedure, called parametric model embedding (PME) \cite{serani2023parametric}, has been recently proposed. The PME simply extends the design-space dimensionality reduction procedure based on KLE/PCA using a generalized feature space that includes shape modification and design variables vectors together with a generalized inner product, building an embedded model of the original design parameterization.

%The analysis of solver trends provides valuable insights into the changing dynamics of computational tool usage in SBDO studies. It reflects both the technological evolution in the field and the shifting focus of optimization challenges over the past few decades, showing how this sector has recently turned towards a deep use of high-fidelity solvers. 
The choice of numerical solvers in SBDO studies reflects an evolving landscape. The growing preference for RANS solvers over potential flow methods marks a shift towards more comprehensive fluid dynamics modeling. This transition aligns with the industry's push towards capturing more complex, turbulent flows and the increasing availability of computational resources. However, the consistent but limited use of FEM solvers indicates a potential underutilization in structural optimization aspects of marine engineering. Future research could benefit from a more integrative approach that combines RANS for fluid dynamics with FEM for structural analysis, potentially leading to more comprehensive and effective optimization in marine engineering.

In the field of engineering optimization, the emphasis is often on achieving an optimal solution in a single iteration of an algorithm, reflecting the practical constraints of time and resources. Traditional stochastic global methods, while robust in exploratory capacity, typically require multiple iterations to ascertain solution reliability due to their inherent randomness. This necessitates a shift towards deterministic variants of global evolutionary strategies and population-based methods. These deterministic adaptations aim to retain the broad exploratory characteristics of global methods but enhance the consistency and predictability of outcomes in each individual run.
Additionally, the strategic integration of these deterministic global methods with deterministic local search techniques marks a significant advancement in optimization practice. This hybrid approach synergistically merges the expansive exploration capabilities of global methods with the focused, efficient refinement of local optimization techniques, such as gradient-based or line search methods. The result is an approach that effectively leverages the strengths of both methodologies, facilitating convergence to the most optimal solution within the constraints of a single algorithmic execution.
Such developments in deterministic global methods, complemented by hybridization with local searches, are particularly salient in engineering contexts. They offer a streamlined and effective means of identifying the global optimum, aligning with the practical exigencies of engineering optimization where timely and reliable solutions are paramount.

The trends and preferences in optimization algorithms and surrogate methods in SBDO reflect an evolving field that continually adapts to the intricacies of marine engineering design problems. The shift towards global optimization and the increasing reliance on surrogate-based methods indicate a strategic response to the challenges of high-dimensional, complex design spaces. This evolution underscores the industry's commitment to finding a balance between computational efficiency and the need for thorough, accurate design exploration. It can be noted how the extension to multi-fidelity approaches, as well as, the integration of active learning/adaptive sampling procedure for the surrogate training process, is still limited. These two branches represent a pathway to follow for future research to assess clearly the pros and cons of multi-fidelity versus single-fidelity methods, as well as identify the most efficient and effective DoE in combination with active learning/adaptive sampling procedure. It may be emphasized that, as for the problem formulation, the literature presents several unclear statements on which DoE is used for surrogate training, as well as what kind of acquisition function has been used in the case of active learning. This represents a huge gap in interpretability and repeatability of the methodologies, that have to be filled. 

Finally, the current distribution of SBDO applications in marine engineering indicates a strong focus on vehicle design, reflecting both the industry's needs and the maturity of optimization techniques in this area. However, the presence of renewable energy and offshore applications, although smaller in proportion, is significant. It suggests a growing recognition of SBDO's potential in these areas, particularly in response to global trends toward sustainable energy and the need for environmentally resilient offshore infrastructure.
As the field of SBDO evolves, it might be expected to see a diversification in its applications. The renewable energy sector, in particular, may experience growth in SBDO applications as the demand for sustainable energy solutions increases. Furthermore, advancements in SBDO methodologies might lead to greater adoption in offshore applications, addressing the unique challenges posed by these environments.
The distribution of SBDO applications across different types of marine vehicles reflects the diverse challenges and priorities in marine vehicle design. The prominence of SBDO in surface vehicle optimization aligns with the global scale and economic significance of these vessels. The focus on underwater vehicles highlights the technological advancements and specialized requirements in this sector. Meanwhile, the application in amphibious vehicle design, although likely less in comparison, underscores the complexity and innovation in multi-environment vehicle design. SBDO is a crucial tool in advancing the design and performance of various types of marine vehicles, addressing unique challenges, and contributing to the evolution of more efficient, capable, and environmentally friendly marine transportation and exploration technologies.
The breakdown of SBDO applications across various types of surface ships demonstrates the versatility and significance of optimization techniques in addressing the diverse design and operational challenges of different ship categories. The focus on containerships and military vessels reflects economic and strategic priorities, while the emphasis on multi-hulls indicates an interest in innovative hull designs. The diverse application across other ship types, such as tankers, bulk carriers, fishing, and passenger ships, highlights the broad applicability of SBDO in enhancing various aspects of marine vessel design and operation.

In summary, while SBDO has become a cornerstone in marine engineering, there is a clear path forward for further advancements. Embracing complex optimization methodologies, expanding the use of MDO, and integrating various computational solvers could pave the way for more innovative and sustainable solutions in marine engineering. These developments, coupled with the broader trends in digitalization and environmental consciousness, are poised to significantly shape the future of SBDO in this field.

{It finally should be noted that while focusing exclusively on peer-reviewed journal papers has ensured the academic rigor and reliability of the sources reviewed, this approach may have limited the representation of industrial applications of SBDO in marine engineering. Industrial projects, especially those involving multi-objective and constrained optimization problems as well as multi-disciplinary efforts, are often not documented in the academic databases surveyed. This is due to various factors, including proprietary considerations and the publication venues typically preferred by industry practitioners, such as industry magazines, conference contributions, and books detailing larger research and development projects.}

%=================================================================================
\section{Conclusions}\label{sec13}

The scoping review conducted in this study underscores the increasingly pivotal role of simulation-based design optimization (SBDO) in marine engineering. Our findings illuminate how SBDO is not just a facilitator of improved performance and cost-efficiency in marine engineering systems and components but also a catalyst for innovation and adaptation in the face of evolving technological and environmental challenges.

Significantly, our analysis reveals a low use of more sophisticated, multi-objective, and stochastic optimization approaches in SBDO, despite the complex, dynamic nature of marine environments. There remains a predominant reliance on simpler, deterministic single-objective formulations, highlighting a crucial area for future development. This gap underscores the necessity for more advanced algorithms that can more accurately model and navigate the uncertainties inherent in marine engineering, including factors like wave dynamics and ocean currents.

Moreover, the review highlights the emergence of high-fidelity solvers in SBDO, signaling a shift towards more nuanced and detailed simulation capabilities. This advancement is indicative of the field's progression towards tackling more complex optimization challenges, further driven by the integration of active learning and adaptive sampling techniques in surrogate-based optimization models and the development of design-space dimensionality reduction procedures for addressing the curse of dimensionality problem.

In conclusion, this scoping review not only reaffirms the significant potential of SBDO in revolutionizing marine engineering practices but also identifies critical pathways for future research. These include the need for more integrative, multidisciplinary approaches, and the development of optimization methods that are both computationally efficient and robust in the face of the unique challenges posed by the marine environment. As the field continues to evolve, these insights will be instrumental in guiding the next generation of research and innovation in SBDO, paving the way for more sustainable, efficient, and advanced marine engineering solutions.

\section*{Acknowledgments}
Dr. Serani has been partially supported by the Horizon Europe 
``RETROFIT55 - Retrofit solutions to achieve 55\% GHG reduction by 2030'', grant agreement 101096068. Dr. Scholcz is grateful to the Dutch Ministry of Economic Affairs which has partially funded the present work.

%Bibliography
\bibliographystyle{unsrt}  
\bibliography{biblio}

\end{document}